\documentclass[a4paper,12pt,reqno,final]{amsart}
\usepackage{float}
\usepackage{graphicx}
\usepackage{amsmath,amsthm,amssymb,amsfonts}
\usepackage{subcaption}
\usepackage{showkeys}  
\usepackage{latexsym}
\usepackage{color}     
\usepackage{hyperref}
\newcommand{\Deg}{\operatorname{deg}}
\newcommand{\interior}{\operatorname{int}}
\newcommand{\imloc}{\operatorname{im}_\loc}
\newcommand{\imT}{\operatorname{im}_T}
\newcommand{\cO}{\mathcal{O}}
\newcommand{\cG}{\mathcal{G}}
\newcommand{\cR}{\mathcal{R}}

\newcommand{\AI}{\operatorname{AI}}
\newcommand{\AIB}{\operatorname{AIB}}
\newcommand{\AIloc}{\AI_\loc}
\newcommand{\AIBloc}{\AIB_\loc}
\newcommand{\CNC}{\operatorname{CNC}}
\newcommand{\INV}{\operatorname{INV}}
\newcommand{\DEGone}{\operatorname{DEG1}}
\newcommand{\DEGoneloc}{\operatorname{DEG1}_\loc}

\usepackage[T1]{fontenc}
\usepackage{amsmath}                       
\usepackage{amssymb}                       
\usepackage{amsfonts}                      
\usepackage{amsthm}     
\usepackage{enumerate}

\usepackage[T1]{fontenc}
\usepackage{amsmath}                       
\usepackage{amssymb}                       
\usepackage{amsfonts}                      
\usepackage{amsthm}     
\usepackage{enumerate}

\newcommand{\RR}{\mathbb R}
\newcommand{\NN}{\mathbb N}
\newcommand{\ZZ}{\mathbb Z}

\newcommand{\cA}{\mathcal{A}}
\newcommand{\cB}{\mathcal{B}}

\newcommand{\cL}{\mathcal{L}}

\newcommand{\cH}{\mathcal{H}}

\newcommand{\cU}{\mathcal{U}}
\newcommand{\eps}{\varepsilon}

\newcommand{\cof}{\operatorname{cof}}
\newcommand{\loc}{\text{loc}}

\newcommand{\mbf}[1]{\boldsymbol{#1}}
\newcommand{\To}{\longrightarrow}

\newcommand{\abs}[1]{\left\vert#1\right\vert}

\newcommand{\norm}[1]{\left\|#1\right\|}
\newcommand{\normb}[1]{\big\|#1\big\|}

\newcommand{\normn}[1]{\|#1\|}
\newcommand{\dist}[2] {\operatorname{dist}\left(#1;#2\right)}
\newcommand{\Dist}[2] {\operatorname{Dist}\left(#1,#2\right)}
\newcommand{\mysetr}[2] {\left\{#1\,\left|\,#2\right.\right\}}
\newcommand{\mysetl}[2] {\left\{\left.#1\,\right|\,#2\right\}}
\newcommand{\mysetb}[2] {\big\{#1\,\big|\,#2\big\}}

\newcommand{\identity}{\operatorname{id}}

\newcommand{\im}{\operatorname{im}}

\newcommand{\be}{\begin{equation}}
\newcommand{\ee}{\end{equation}}
\newcommand{\bald}{\begin{aligned}}
\newcommand{\eald}{\end{aligned}}
\newcommand{\baldat}{\begin{alignedat}}
\newcommand{\ealdat}{\end{alignedat}}

\newenvironment{myenum}{\vspace*{-.5ex}\begin{enumerate}}{\end{enumerate}\vspace*{-.5ex}}

\newtheorem{theorem}{Theorem}%[]

\newtheorem{thm}[theorem]{Theorem}
\newtheorem{cor}[theorem]{Corollary}
\newtheorem{lem}[theorem]{Lemma}
\newtheorem{prop}[theorem]{Proposition}
\theoremstyle{definition}
\newtheorem{defn}[theorem]{Definition}
\newtheorem{ex}[theorem]{Example}
\theoremstyle{remark}
\newtheorem{rem}[theorem]{Remark}
\newtheorem{problem}[theorem]{Problem}
\renewcommand{\proofname}{\bfseries{Proof}}
\renewenvironment{proof}[1][\proofname]{\par
  \normalfont
  \trivlist
  \item[\hskip\labelsep\itshape
    \bfseries{#1.}]\ignorespaces
}{%
  \qed\endtrivlist
}

\begin{document}

\title[]{Global invertibility for orientation-preserving Sobolev maps via invertibility on or near the boundary}
\
\author{Stefan Kr\"omer} \address{Stefan Kr\"omer, The Czech Academy of Sciences, Institute of Information Theory and Automation, Pod vod\'{a}renskou v\v{e}\v{z}\'{\i}~4, 182~08~Praha~8, Czech Republic, \email{}{skroemer@utia.cas.cz}}
\date{\today}

\begin{abstract} 
By a result of John Ball (1981), a locally orientation preserving Sobolev map is almost everywhere globally invertible whenever its boundary values admit a homeomorphic extension. As shown here for any dimension, the conclusions of Ball's theorem and related results can be reached while completely avoiding the problem of homeomorphic extension. For suitable domains, it is enough to know that the trace is invertible on the boundary or can be uniformly approximated by such maps.
An application in Nonlinear Elasticity is the existence of homeomorphic minimizers with finite distortion whose boundary values are not fixed. As a tool in the proofs, strictly orientation-preserving maps and their global invertibility properties are studied from a purely topological point of view. 
\end{abstract}

\maketitle


\maketitle

%*************************************************\\
%TODO:
%END TODO\\
%*************************************************\\

\section{Introduction\label{sec:intro}}

A classical problem in Nonlinear Elasticity is to determine
whether a Sobolev map $y:\Omega\to \RR^d$ on a bounded domain $\Omega\subset \RR^d$, 
is invertible in a suitable sense. In this context, the map $y$ describes the deformation of 
an elastic solid occupying $\Omega$ in its undeformed state. In this model,
lack of invertibilty corresponds to self-interpenetration which is clearly undesirable. If we assume the existence of a stored energy density (i.e., elastic deformation does not dissipate energy) with suitable properties,
then a stable deformed state can be found via global minimization \cite{Ba77a}.
The typical examples for such energies enforce $y\in W^{1,p}(\Omega;\RR^d)$, strictly orientation-preserving in the sense that $\det \nabla y>0$ a.e.~in $\Omega$.

However, even if $y\in C^1$ and $\det \nabla y$ is positive everywhere, this does not suffice to guarantee global invertibility, because different ends of the body can still overlap.
The variational theory is compatible with imposing global invertibility (in a weak a.e.~sense) as a constraint, the Ciarlet-Ne\v{c}as condition \eqref{CN0} \cite{CiaNe87a} (see Def.~\ref{def:CNC}).

In case of a strictly orientation-preserving map, if we also assume that its boundary values match those of a homeomorphism, \eqref{CN0} always holds as a byproduct of a result of Ball \cite{Ba81a}. 
The result also provides further topological properties of the deformation and its image. 
For similar purposes, this assumption also appears in \cite{Sve88a} and other works, in particular in context of maps of finite distortion \cite{HeKo14B} (e.g.).
The caveat here is that for $d \geq 3$, a homeomorphic extension can fail to exist. This issue is discussed in detail in Appendix~\ref{sec:A:homext}.

We will see that for the purpose of proving the result of \cite{Ba81a},
it suffices to know that the deformation is continuous and 
invertible on the boundary, 
or, more generally, approximately invertible on the boundary ($\AIB$) with respect to uniform convergence (Theorem~\ref{cor:Ball}). 
The result requires an assumption on the topological nature of the domain, namely, that $\RR^d\setminus \partial\Omega$ consists of exactly two connected components. 
As illustrated by examples in Appendix~\ref{sec:A:topological}, this restriction is not just a technical issue.
The class $\AIB$ can also be used instead of \eqref{CN0} to implement a global invertibility constraint in Nonlinear Elasticity, because it is stable under weak convergence in suitable Sobolev spaces and therefore compatible with direct methods in the Calculus of Variations (Section~\ref{sec:AI}).

A first, crucial step to connect global invertibility with invertibility on the boundary amounts to calculating the topological degree of $y$  (Theorem~\ref{thm:deg-of-def} in Section~\ref{sec:AIBdeg}).
It is of course not surprising that information on the boundary suffices for that, since the degree only depends on boundary values. Actually,
if we assume in addition that $y(\bar\Omega)=\bar\Omega$ and $y(\partial\Omega)=\partial\Omega$, then 
we are working in a class of maps with a semigroup structure, and it is well known that homeomorphisms in such a class 
have a degree of $\pm 1$, due to a standard multiplicativity property of the degree. 
The main global assumption of \cite{Ba81a}, the existence of a homeomorphic extension,
would allow us to reduce our more general situation to this scenario. 
%As a matter of fact, 
%this is even the only point where the existence of the homeomorphic extension is truly exploited in \cite{Ba81a}. 

However, in general there is no natural group structure we could use directly. 
Even if there exists a degree for endomorphisms of $y(\partial\Omega)$ (neither $\partial\Omega$ nor $y(\partial\Omega)$ 
are always topological $(d-1)$-manifolds in our setting!), we still have the problem of linking it
with the definition of the degree for continuous maps $\RR^d\to\RR^d$. This is not trivial, as any comparison of the two and their necessary normalizing conditions 
(like $\Deg(\identity;\Omega;\cdot)=1$ in $\Omega$) already requires a continuous map $\partial\Omega\to y(\partial\Omega)$ with 
a continuous extension such that \emph{both degrees are known}, to act as reference map which is meaningful in both worlds.
Instead, the proof of Theorem~\ref{thm:deg-of-def} exclusively works with the degree for continuous endomorphisms of $\RR^d$ on domains. It relies on a generalized version of the Jordan-Brouwer separation theorem and a general formula for the degree of composite functions, the multiplication theorem (see \cite{FoGa95B}, e.g.). 

As a second tool for the proof of Theorem~\ref{cor:Ball} and further applications concerning global invertibility in $W^{1,p}$ for $p\geq d$ (Section~\ref{sec:p>d}), we develop a self-contained, purely topological theory of strictly orientation preserving maps with the help of Brouwer's degree (Section~\ref{sec:top}). 
Their interplay with a global invertibility constraint stated in terms of the degree is summarized in Theorem~\ref{thm:sop-deg1} and Corollary~\ref{cor:sop-deg1loc}. 
Besides the generalization of \cite{Ba81a}, this can be used to complement results of
\cite{ViMa98a} and \cite{Raj10a,Raj11a} for maps of finite distortion to 
prove that suitable deformations are actually homeomorphisms (Theorem~\ref{thm:finoutdistort} and Theorem~\ref{thm:fininnerdistort}). As a direct application, we show the existence of homeomorphic minimizers for nonlinear elastic energy functionals controlling the inner or outer distortion (Subsection~\ref{subsec:min}), without fixing the boundary values of admissible deformations (in the spirit of \cite{Ba81a}, as in \cite{HeKo14B}, e.g.) or prescribing a given Lipschitz domain as their image (as in \cite{MoVo19a}). 

Another application will be presented in the forthcoming paper \cite{KroeVa20P}, in context of the numerical approximation 
of global invertibility constraints via penalty terms in the energy. One such approach is given in \cite{KroeVa19a}, where the penalty term acts on the full domain. However, reflecting the nonlocal nature of global invertibility,
any suitable penalty term is necessarily nonlocal, and the associated computational cost can be significantly reduced by using a variant acting only on the boundary \cite{KroeVa20P}, thereby reducing the effective dimension of the problem. A priori, such a boundary penalty term can only hope to ensure invertibility on the boundary, which is why we need the results developed here to understand the link to full invertibility.

%From a technical point of view, the results of Section~\ref{sec:top} identify and unify the local %and global topological aspects of strictly orientation preserving maps, separating them from %additional analytic fine properties of functions in the respective investigated class of several %approaches in the literature, in particular  \cite{Sve88a}, \cite{MueSpe95a}, \cite{HeMo10a} and %\cite{BaHeMo17a}. A generalization capable of universally handling various examples of suitable %classes in $W^{1,p}$ with $d-1<p<d$ as treated in these articles might be possible, ideally %clarifying their topological interconnections and allowing for cavities. This is not trivial, %though, and can hopefully be addressed in future work.

All results presented here apply in particular for $d=3$.

\tableofcontents

\subsection{Basic notation and terminology}

Throughout the article, we use $p\geq d\geq 2$, where $p\in \RR$ and $d\in \NN$, and the 
following subsets of the Sobolev space $W^{1,p}(\Omega;\RR^d)$
of functions on an open set $\Omega\subset \RR^d$ with values in $\RR^d$:
\begin{align*}
\baldat{2}
	&W^{1,p}_+(\Omega;\RR^d)&&:=\mysetr{y\in W^{1,p}(\Omega;\RR^d)}{\text{$\det\nabla y>0$ a.e.~in $\Omega$}}, \\
	&W^{1,p}_{+,\loc}(\Omega;\RR^d)&&:=\mysetr{y\in W^{1,p}_\loc(\Omega;\RR^d)}{\text{$\det\nabla y>0$ a.e.~in $\Omega$}}.
\ealdat
\end{align*}
Here, $\Omega$ is called a \emph{domain} if is open and connected, and "a.e." abbreviates almost everywhere, which is understood with respect
to the Lebesgue measure $\cL^d$ in $\RR^d$ unless specified otherwise.
For any set $A\subset\RR^d$, $\bar{A}$ is its closure and $\interior A$ its interior, and for $A_1,A_2\subset \RR^d$,
$A_1\subset\subset A_2$ means that $\bar{A}_1$ is compact and $\bar{A}_1\subset A_2$.
An Euclidean norm is always denoted by $\abs{\cdot}$, in any finite-dimensional real vector space that should be clear from the context,
and for a point $x$ and a set $S$ in such a space, $\dist{x}{S}:=\inf \mysetl{\abs{x-s}}{s\in S}$.

Notation concerning the topological degree is introduced in Section~\ref{sec:degreebasic}.

\section{Constraints related to global invertibility} \label{sec:AI}

In this section, we collect various conditions related to global invertibility that are viable as constraints for variational approaches in Nonlinear Elasticity. It is interesting to note that even in the "simple" case of $p>d$ on smooth domains,
it is not clear whether the notions based on approximate invertibility 
coincide with or are stronger than the more classical constraints like the Ciarlet-Ne\v{c}as condition, cf.~Remark~\ref{rem:invcompare}. 

\subsection{Approximate invertibility}

\begin{defn}[$\AI$: approximately invertible on a compact set]\label{def:AI}~\\
Let $K\subset \RR^d$ be bounded. A continuous function $y:K\to \RR^d$ is called \emph{approximately invertible on $K$} if there exists a sequence of injective maps
$\varphi_k\in C(K;\RR^d)$ with $\varphi_k\to y$ uniformly on $K$.
The class of all such maps $y$ in $C(K;\RR^d)$ is denoted by $\AI(K)$. 
\end{defn}
The most important examples are $K=\bar\Omega$ and $K=\partial\Omega$ on a domain $\Omega\subset\RR^d$.
As a matter of fact, the class $\AI(K)$ is linked to monotone mappings in the topological sense. 
For these, in dimension $d=2$ (but not for $d\geq 3$), a quite comprehensive and satisfactory theory is available \cite{IwaOnn16a,IwaOnn19a}.
In this article, we are especially interested in $d\geq 3$ and $K=\partial\Omega$:
\begin{defn}[$\AIB$: approximately invertible on the boundary]\label{def:AIB}~\\
Let $\Omega\subset \RR^d$ be open and bounded, and $y:\bar\Omega \to \RR^d$ with $y\in C(\partial\Omega;\RR^d)$. 
We say that $y$ is \emph{approximately invertible on the boundary}, if 
$y\in \AIB:=\AI(\partial\Omega)$. 
\end{defn}
From the point of view of Nonlinear Elasticity, $\AIB$ describes the class of deformations whose
deformed boundaries can be moved out of self-contact.

It is easy to see that $\AI$ is sequentially closed under uniform convergence on $K$.
In suitable Sobolev spaces, this implies stability under weak convergence:
\begin{lem}\label{lem:AIPwclosed}
Let $y\in W^{1,p}(\Omega;\RR^d)$ and let $(y_k)\subset W^{1,p}(\Omega;\RR^d)$ be a sequence
with $y_k\rightharpoonup y$ weakly in $W^{1,p}$.
In addition, let $K\subset \bar\Omega$ be compact and assume that $(y_k)\subset \cA\cap \AI(K)$ for a set $\cA\subset W^{1,p}(\Omega;\RR^d)$ 
compactly embedded in $C(K;\RR^d)$.
Then $y\in \AI(K)$.
\end{lem}
\begin{proof}
This follows from a straightforward diagonal subsequence argument in $C(K;\RR^d)$.
\end{proof}
\begin{rem}\label{rem:AIwlsc}
The assumptions of Lemma~\ref{lem:AIPwclosed} are satisfied in each of the following cases:
\begin{enumerate}
\item[(a)] $p>d$, $K\subset\Omega$ and $\cA=W^{1,p}$;
\item[(b)] $p>d$, $\Omega$ is a Lipschitz domain, $K\subset\bar{\Omega}$ and $\cA=W^{1,p}$;
\item[(c)] $p\geq d$, $K\subset \Omega$ and $\cA=W_+^{1,p}=W^{1,p}\cap \{\det \nabla y>0~\text{a.e.}\}$.
\item[(d)] $p>d-1$, $\cA=\{y\}\cup \mysetr{y_k}{k\in\NN}$ and $K\subset \Omega$ s.t.~
$\cA$ is embedded and closed in $C(K;\RR^d)$.
\end{enumerate}
Here, (a) and (b) are due to standard compact embeddings of $W^{1,p}$. For (c) see Remark~\ref{rem:pequald}.
In (d), one has to be careful to correctly interpret $\cA$ as a subset of $C(K;\RR^d)$ in a way independent of the choice of representatives. 
For instance, if we fix $x_0\in\Omega$, then $K_r:=\partial B_r(x_0)$ is admissible in (d) for $\cL^1$-a.e.~$r\in (0,\dist{x_0}{\partial\Omega})$,
essentially because we work with countably many functions and for a.e.~$r$, $y_k\to y$ weakly in $W^{1,p}(K_r;\RR^d;\cH^{d-1})$ (the Sobolev space with respect to the surface measure $\cH^{d-1}$). For more details see, e.g., \cite{MueSpe95a}.
\end{rem}
As a consequence, we can easily obtain the existence of minimizers in $\AI(K)$,
for instance for 
the nonlinear elastic energies with polyconvex energy density studied by Ball \cite{Ba77a} and M\"uller \cite{Mue90a}:
\newcommand{\cY}{\mathcal{Y}}
\begin{thm}\label{thm:min}
Let $p\geq d$, let $\Omega\subset \RR^d$ be a Lipschitz domain, let 
$E: W^{1,p}(\Omega;\RR^d)\to \RR\cup\{+\infty\}$ such that $E(y)=+\infty$ for all $y\notin W^{1,p}_+$, i.e., whenever $\det\nabla y\leq 0$ on a set of positive measure. In addition, assume that
$E$ is coercive and weakly sequentially lower semincontinuous in $W^{1,p}$.
If $E\not\equiv +\infty$ and $K\subset \Omega$ is compact, then $E$ attains its minimum on $\cY:=\AI(K)\cap W^{1,p}_+(\Omega;\RR^d)$.
Moreover, if $p>d$, the above also holds for compact $K\subset \bar\Omega$.
\end{thm}
\begin{proof}
This follows by the direct method.
%This immediately follows by the Direct Method: 
%For any $y_n\in \cY$ with $E(y_n)\to I:=\inf_{y\in \cY} E(y)$,
%$(y_n)$ is bounded in $W^{1,p}$ by coercivity of $E$. Hence, up to a subsequence, $y_n\rightharpoonup y^*$ in $W^{1,p}$,
%and by weak lower semicontinuity, $E(y^*)\leq I$. Finally, by Lemma~\ref{lem:AIPwclosed}, $y^*\in \AI(K;\RR^d)$. Hence, $y^*\in \cY$ is admissible and a minimizer.
\end{proof}
\begin{rem}
Here and in the rest of the article, it is implicitly understood 
that we always use the continuous representative of $y$ if available.
\end{rem}

All functions which admit homeomorphic extensions from the boundary into the domain are in $\AIB$:
\begin{prop}\label{prop:homextAIB}
Let $\Omega\subset \RR^d$ be a bounded Lipschitz domain. 
If $y:\bar\Omega\to \RR^d$ is continuous, and $y|_{\Omega}:\Omega\to y(\Omega)$ is invertible, then
$y\in \AIB$.
\end{prop}
\begin{proof}
Since $\Omega$ is a Lipschitz domain, there exists a sequence of invertible maps $\Psi_k:\bar\Omega\to \Psi_k(\bar\Omega)\subset \RR^d$ of class $C^\infty$ 
such that $\Psi_k(\Omega)\subset\subset\Omega$ (slightly smaller), while $\Psi_k\to \identity$ in $W^{1,\infty}$ as $k\to\infty$.
Locally, in each cube where the boundary is represented as the graph of a Lipschitz function, $\Psi_k$ can be defined as the affine map slightly shrinking the local piece of $\Omega$ 
"down" into itself, and since all these maps are still close to the identity, they can be easily glued by a smooth decomposition of unity.
Thus, $\varphi_k:=y\circ \Psi_k|_{\partial\Omega}$ is a sequence of continuous, injective maps with $\varphi_k\to y|_{\partial\Omega}$ in $C(\partial\Omega;\RR^d)$.
\end{proof}
\begin{rem}
The converse of Proposition~\ref{prop:homextAIB} includes the problem of homeomorphic extension as a special case, namely, if we only consider $y\in \AIB$ which is already invertible on $\partial\Omega$. In general, it is false for $d\geq 3$ because not all such $y$ admit 
a homeomorphic extension into the domain, not even if $y|_{\partial\Omega}$ is bi-Lipschitz (see Remark~\ref{rem:Fox-Artin}).
\end{rem}
We can also slightly modifiy the definition of $\AIB$, allowing the approximants to be defined on sets approaching $\partial\Omega$ from the
inside. Monotone coverings of $\Omega$ from the inside with a mild regularity property (to be used later) are helpful for that purpose:
\begin{defn}[regular inner covering] \label{def:regcover}
Let $\Omega\subset \RR^d$ be open. We call a family of sets $(\Omega_m)_{m\in\NN}$ a \emph{regular inner covering of $\Omega$},
if $\Omega_m$ is open and bounded, $\Omega_m\subset \Omega_{m+1} \subset\subset \Omega$, $\bigcup_{m\in\NN}\Omega_m=\Omega$ and 
$\cL^d(\partial\Omega_m)=0$ for every $m\in\NN$.
\end{defn}
\begin{rem}\label{rem:regcover}
It is not difficult to see that a regular inner covering always exists; as a matter of fact, we could even assume that $\partial\Omega_m$ is of class $C^\infty$ instead of just being a set of measure zero. Moreover, if we know that $\RR^d\setminus \partial\Omega$ has only two connected components (which is important to apply Theorem~\ref{thm:deg-of-def} below), 
we can always find a regular inner covering such that all $\Omega_m$ inherit this property.
\end{rem}
The following variants of $\AIB$ and $\AI$ are particularly 
useful when $y$ is continuous in $\Omega$ but cannot be continuously extended to $\partial\Omega$ 
(like maps in $W_+^{1,p}$ with $p=d$):
\begin{defn}[$\AIBloc$, $\AIloc(\Omega)$]\label{def:AIBloc}
Let $\Omega\subset \RR^d$ be open, let $y:\Omega\to \RR^d$ be continuous and let 
$(\Omega_m)_{m\in\NN}$ be regular inner covering of $\Omega$.
We say that $y\in \AIBloc$ with respect to $(\Omega_m)$,
if for each $m\in\NN$, $y\in \AIB$ on $\Omega_m$, i.e.,
if there exists continuous and injective maps $\varphi^{(m)}_k:\partial\Omega_m\to \RR^d$ such that
\[
	\normb{y-\varphi^{(m)}_k}_{C(\partial\Omega_m)}\to 0\quad\text{as $k\to\infty$}.
\]
Analogously, we say that $y\in \AIloc(\Omega)$ with respect to $(\Omega_m)$, if for each $m\in\NN$, $y\in \AI(\Omega_m)$.
\end{defn}
\begin{rem}\label{rem:AIBloc}
It is easy to see that both Lemma~\ref{lem:AIPwclosed} and Theorem~\ref{thm:min} also hold for
$\AIBloc$ or $\AIloc(\Omega)$ instead of $\AI(K)$. In fact, we do not even use that $(\Omega_m)$ is a countable family, because there are no conditions
linking $\varphi^{(m)}_k$ for two different values of $m$. 
This means that diagonal subsequences as in the proof of Lemma~\ref{lem:AIPwclosed} can be chosen for each fixed $m$ separately, by the axiom of choice if we have more than countably many $m$.
\end{rem}
\begin{rem}
$\AIB_\loc$ can potentially still be used in
settings with low regularity like $W^{1,p}$ with $d-1<p<d$, cf.~Remark~\ref{rem:AIwlsc}.
By contrast, this does not work for $\AI_\loc(\Omega)$ as defined here. 
However, approximate invertibility defined with respect to weak convergence in the Sobolev space is a still meaningful concept \cite{BouHeMa19Pa}.
\end{rem}
An explicit example for an existence results in the spirit of Theorem~\ref{thm:min} using either $\AIBloc$ or $\AIloc(\Omega)$ as a constraint is given in Subsection~\ref{subsec:min}.

Unlike all the other invertibility constraints presented in this section, $\AIB$ and $\AIBloc$ only restrict the
global behavior of $y$ with information given on or near the boundary, but not its local properties inside the domain. 
However, in $W_+^{1,p}$ with $p\geq d$, local restrictions still follow automatically, see Remark~\ref{rem:invcompare}. Practically, $\AIB$ can be easier to show for a given deformation, though.

\subsection{The Ciarlet-Ne\v{c}as condition and condition (INV)}

For comparison, we briefly recall two invertibility constraints often used in the literature. 

The standard constraint for $p>d$ was introduced in \cite{CiaNe87a}:
\begin{defn}[Ciarlet-Ne\v{c}as condition \eqref{CN0}]\label{def:CNC}~\\
Let  $\Omega\subset \RR^d$ be open and bounded. 
A map $y\in W^{1,p}_+(\Omega;\RR^d)$, $p\geq d$, satisfies the \emph{Ciarlet-Ne\v{c}as condition},
or, shortly, $y\in \CNC$, if
\begin{align*}\label{CN0}\tag{CNc}
	\int_\Omega \det\nabla y(x)\,dx\leq \cL^d(y(\Omega)).
\end{align*}
\end{defn}
Using the area formula as in the proof of Lemma~\ref{lem:CNloc} below, it is not difficult to see that
\eqref{CN0} is equivalent to injectivity almost everywhere in the sense that
the set of all points in $\Omega$ where $y$ fails to be injective is of measure zero.

If $p>d$ and $\Omega$ is Lipschitz, \eqref{CN0} is stable under weak convergence in $W^{1,p}$ for $p>d$ (proved as part of \cite[Theorem 5]{CiaNe87a}).
If $p=d$, it can happen that the left hand side of \eqref{CN0} jumps down in the limit along a weakly converging sequence, due to concentration 
effects at the boundary (cf.~\cite{Mue90a}, \cite{KaKroeKru14a}). Nevertheless, the right hand side actually produces a matching jump in such cases because
\eqref{CN0} still behaves stably as a whole along sequences in $W^{1,d}_+$
weakly converging in $W^{1,d}$, a result obtained in \cite{Ta88a} in broader context. 
Alternatively, one can use the even more general results of \cite{GiaPo08a},
or Remark~\ref{rem:invcompare} (e) below.

The following lemma also used in \cite{Ta88a} shows that there is no point in defining a "loc" version of \eqref{CN0}. 
In view of all the properties known for functions in $W^{1,d}_+$ in the interior of the domain (cf.~Remark~\ref{rem:pequald}),
and the local equi-integrability result for the determinant of \cite{Mue90a}, it
can also be the basis of yet another, more direct proof of the result of \cite{Ta88a} in the special case $W^{1,d}_+$.
\begin{lem}\label{lem:CNloc}
Let $\Omega\subset \RR^d$ be open and bounded, let $y\in W^{1,d}_+(\Omega;\RR^d)$
and let $(\Omega_m)_{m\in\NN}\subset\Omega$ be a sequence of open sets with $\Omega_{m}\subset \Omega_{m+1}\subset\subset \Omega$
and $\textstyle \bigcup_{m\in\NN} \Omega_{m}=\Omega$. Then $y$ satisfies \eqref{CN0} on $\Omega$ if and only if it satisfies \eqref{CN0} on $\Omega_m$ for all $m$.
\end{lem}
\begin{proof}
By the area formula (see, e.g., \cite[Theorem 5.34]{FoGa95B}),
\begin{align} \label{areaformula0}
	\int_\Omega \abs{\det\nabla y}\,dx = \int_{y(\Omega)} \#y^{-1}(\{z\})\,dz,
\end{align}
where $\#y^{-1}(\{z\})$ denotes the number of elements of $y^{-1}(\{z\})$. 
In view of \eqref{areaformula0}, \eqref{CN0} can be expressed as $\#y^{-1}(\{z\})=1$ for a.e.~$z$,
and it is clear that \eqref{CN0} on $\Omega$ implies \eqref{CN0} on every smaller set like $\Omega_m$.
The converse follows by monotone convergence.
\end{proof}

Another condition implying both local and global invertibility properties 
was developed in \cite{MueSpe95a}, mainly intended for settings involving cavitation
in $W^{1,p}$ with $d-1<p<d$. It uses the fact that the degree only depends on boundary values, cf.~Section~\ref{sec:degreebasic}, and
the concept of the topological image $\imT$ based on the degree, cf.~Lemma~\ref{lem:orpres-top}.
\begin{defn}[M\"uller-Spector condition \eqref{INV}]\label{def:INV}
Let $\Omega\subset \RR^d$ be a bounded domain and $y:\Omega\to \RR^d$. 
%such that $y|_{\partial B_r(a)}$ is continuous for all $a\in\Omega$ and a.e.~$r\in (0,\dist{a}{\partial\Omega})$.
The map $y$ satisfies \emph{condition \eqref{INV}},
or, shortly, $y\in \INV$, if the following holds:
\begin{align*}\label{INV}\tag{INV}
\bald
	&\text{For every $a\in\Omega$, there exists a set $N_a\subset \RR$ with $\cL^d(N_a)=0$ }\\
	&\text{such that $y\in C(\partial B_r(a);\RR^d)$ for all $r\in (0,\dist{a}{\partial\Omega})\setminus N_a$,}\\
&\baldat[t]{3}
	&\text{(i)}~~&& y(x)\in \imT(y;B_r(a))\cup y(\partial B_r(a))~~\text{for a.e. $x\in \bar B_r(a)$, and}\\
	&\text{(ii)}~~&& y(x)\in \RR^d\setminus \imT(y;B_r(a))~~\text{for a.e. $x\in \Omega\setminus B_r(a)$}.
\ealdat
\eald
\end{align*}
\end{defn}
Condition \eqref{INV} is stable under weak convergence in $W^{1,p}$ for $p>d-1$ \cite[Lemma 3.3]{MueSpe95a}.

\subsection{Maps of topological degree at most one}

Another property that prevents global self-interpenetration in suitable settings can be expressed with the help of the topological degree 
(cf.~Section~\ref{sec:degreebasic}).
This turns out to be a natural common denominator of all the other global invertibility conditions, at least within $W^{1,p}_+$ for $p\geq d$. 
%For the notation and more details including a link to $\AIB$ see Sections~\ref{sec:AIBdeg} and~\ref{sec:top}.
\begin{defn}[Maps of degree at most one]\label{def:deg1}~\\ %\eqref{deg1} and \eqref{deg1b}
Let $\Omega\subset \RR^d$ be open and bounded and $y:\bar\Omega\to \RR^d$.
If $y\in C(\partial\Omega;\RR^d)$, we say that \emph{$y$ is of degree (at most) one}, or, shortly, $y\in \DEGone$, if
\begin{align}\label{deg1}%\tag{deg-1}
	\Deg(y;\Omega;z)\leq 1\quad\text{for all $z\in \RR^d\setminus y(\partial\Omega)$}
\end{align}
If $(\Omega_k)_{k\in\NN}$ is a regular inner covering of $\Omega$ (see Definition~\ref{def:regcover}),
we say that \emph{$y$ is of degree (at most) one locally near the boundary}, or, shortly, $y\in \DEGoneloc$, if
\begin{align}\label{deg1b}%\tag{deg-1b}
	\Deg(y;\Omega_k;z)\leq 1\quad\text{for each $k\in\NN$ and all $z\in \RR^d\setminus y(\partial\Omega_k)$}
\end{align}
\end{defn}
\begin{rem}
As the degree is continuous with respect to uniform convergence on the boundary (stability), 
Lemma~\ref{lem:AIPwclosed} and Theorem~\ref{thm:min} 
also hold if we replace $\AI(K)$ 
by 
\[
\baldat{3}
&\text{(a)}~~ && \DEGone,~~~&&\text{if $p>d$ and $\Omega$ is Lipschitz, or}\\
&\text{(b)}~~ && \DEGoneloc,~~~&&\text{if $p\geq d$.}
\ealdat
\]
In other words, $\DEGone$ and $\DEGoneloc$, too, are stable under weak convergence and viable as variational constraints.
\end{rem}
\begin{rem}\label{rem:invcompare}
Let $\Omega\subset \RR^d$ be a bounded domain. 
Whenever $\AIB$ or $\AIBloc$ are involved, we also assume that $\RR^d\setminus \partial\Omega$ has only two connected components,
to be able to apply Theorem~\ref{thm:deg-of-def}.
Consider the following two classes of strictly orientation preserving Sobolev maps:
\[
\baldat{2}
	&\cY_+^p(\Omega)&&:=W^{1,p}_+(\Omega;\RR^d)\cap C(\Omega;\RR^d)\quad\text{and}\\
  &\cY_+^p(\bar\Omega)&&:=W^{1,p}_+(\Omega;\RR^d)\cap C(\bar\Omega;\RR^d)\cap \mysetr{y}{\cL^d(y(\partial\Omega))=0}
\ealdat
\]
(If $p>d$ and $\Omega$ is a Lipschitz domain, $\cY_+^p(\Omega)=\cY_+^p(\bar\Omega)=W^{1,p}_+(\Omega;\RR^d)$.)
Within $\cY_+^p(\Omega)$ or $\cY_+^p(\bar\Omega)$, respectively, the invertibility conditions are related as follows for $p\geq d$:
\[
\baldat{3}
	&\text{(a)}\quad&& \cY_+^p(\bar{\Omega})\cap \AIB &&\subset~ \cY_+^p(\bar{\Omega})\cap\DEGone;\\
	&\text{(b)}\quad&& \cY_+^p(\bar{\Omega})\cap \CNC&&=~ \cY_+^p(\bar{\Omega})\cap\DEGone;\\
	&\text{(c)}\quad&& \cY_+^p(\bar{\Omega})\cap \DEGoneloc &&=~ \cY_+^p(\bar{\Omega})\cap\DEGone;\\
	&\text{(d)}\quad&& \cY_+^p(\Omega)\cap \AIBloc &&\subset~ \cY_+^p(\Omega)\cap\DEGoneloc;\\
	&\text{(e)}\quad&& \cY_+^p(\Omega)\cap \AIloc(\Omega) &&\subset~ \cY_+^p(\Omega)\cap\DEGoneloc;\\
	&\text{(f)}\quad&& \cY_+^p(\Omega)\cap \CNC&&=~ \cY_+^p(\Omega)\cap\DEGoneloc;\\
	&\text{(g)}\quad&& \cY_+^p(\Omega)\cap \INV&&=~ \cY_+^p(\Omega)\cap\DEGoneloc.
\ealdat
\]
For a proof of some of these connections, we occasionally need properties of the degree (Section~\ref{sec:degreebasic}) and other results presented later. Throughout, Lemma~\ref{lem:orpres-W1p} ensures that any $y\in W^{1,p}_+$ is strictly orientation preserving in the topological sense of 
Section~\ref{sec:top}; in particular, its degree can never be negative.

The inclusions (a) and (d) are consequences of Theorem~\ref{thm:deg-of-def} applied on $\Omega_m$ (see also Remark~\ref{rem:regcover}), and (e) 
analogously follows from Theorem~\ref{thm:deg-of-def2}. In case of (c), "$\subset$" is a consequence of the continuity of the degree (stability) while "$\supset$" follows from 
Corollary~\ref{cor:sop-deg1loc} (i).
Since we also have Lusin's property (N) (cf.~Remark~\ref{rem:pequald}),
$y(\partial \Omega_m)$ has empty interior for all the sets of the regular inner covering of $\Omega$ associated to $\DEGoneloc$.
For (f), one can use Lemma~\ref{lem:CNloc} and the
change-of-variables formula involving the degree \eqref{degCOV} with $f=1$. 
To see (b), we combine (c) and (f) with Lemma~\ref{lem:CNloc}.
Finally, (g) is the content of Lemma~\ref{lem:INVvsDEG} below. 

In (a), (d) and (e), I do not know if the reverse inclusions hold (for $d\geq 3$).
Similarly, while we trivially have that $\AI(\bar\Omega)\subset \AIB$ and $\AIloc(\Omega) \subset \AIBloc$ (with the covering $(\Omega_m)$ fixed), it is not clear if equality holds (given that $\RR^d\setminus \partial\Omega$ has only two connected components). This is related to a weaker variant of the problem of homeomorphic extension for which the counterexamples mentioned in Appendix~\ref{sec:A:homext} do not apply.
\end{rem}

\begin{lem}\label{lem:INVvsDEG}
Let $\Omega\subset \RR^d$ be a bounded domain and $p\geq d$. Then 
\[
	W^{1,p}_+(\Omega;\RR^d)\cap C(\Omega;\RR^d) \cap \INV= W^{1,p}_+(\Omega;\RR^d)\cap C(\Omega;\RR^d)\cap\DEGoneloc.
\]
\end{lem}
\begin{proof}
Let $y\in W^{1,p}_+(\Omega;\RR^d)\cap C(\Omega;\RR^d)$. 
By Lemma~\ref{lem:orpres-W1p}, $y$ is strictly orientation preserving in the topological sense, and by Remark~\ref{rem:pequald},
it satisfies Lusin's property (N). In particular, $y(\partial \Omega_m)$ has measure zero and thus empty interior.

"$\mbf{\subset}$": Suppose that $y$ also satisfies $\eqref{INV}$. 
By \cite[Lemma 3.4]{MueSpe95a}, $y$ is invertible almost everywhere, which implies $y\in\CNC$ by the area formula 
(as in the proof of Lemma~\ref{lem:CNloc}). By Remark~~\ref{rem:invcompare} (f), we conclude that $y\in \DEGoneloc$.

"$\mbf{\supset}$": Given $y\in \DEGoneloc$, Part (i) of \eqref{INV} 
follows from Lemma~\ref{lem:orpres-top} (iii). Part (ii) follows from the fact that by 
Theorem~\ref{thm:sop-deg1} (ii) and (iii) (applied with $U=\Omega_m$, for all $m$), 
$y^{-1}(\{z\})$ can only have more than one connected component if all of its 
components touch $\partial\Omega$ (see also the Remarks~\ref{rem:Ry} and~\ref{rem:Ry2}). Here, $z\in\RR^d$ is arbitrary.
\end{proof}

\section{The degree: basic notation and properties}\label{sec:degreebasic}

In the next two sections, we will heavily use the topological degree (Brouwer's degree). 
We therefore briefly recall its main features. For a definition and its properties see for instance \cite{Kie12B}, \cite{FoGa95B} or \cite{OutRui09B}.

The degree for functions in $\RR^d$ is a number
\begin{align}\label{degreebasic}
	\Deg(y;A;z)\in \ZZ\quad\text{if $z\notin y(\partial A)$},
\end{align}
defined for any continuous map $y:\bar{A}\to \RR^d$ on an open and bounded set $A\subset\RR^d$, with respect to a value $z\in \RR^d$. 
By the Tietze extension theorem, we can always assume that $y:\RR^d\to \RR^d$ is continuous. 
The restriction on the admissible points $z$ in \eqref{degreebasic} is necessary for its definition. Throughout,
it is always assumed to be present, even if not stated explicitly in shorthand notations like $\deg(y;A;\cdot)$.

Besides being integer-valued, the key properties of the degree are the following:
\[
\baldat{2}
&\text{(nomalization)}~~~&& \Deg(\identity;A;z)=1~~\text{for all $z\in A$.} \\
&\text{(additivity)}~~~&& 
\bald[t] 
	&\Deg(y;A_1\cup A_2;z)=\Deg(y;A_1;z)+\Deg(y;A_2;z)\\
	&\quad \text{if $A_1\cap A_2=\emptyset$ and $z\notin \partial A_1\cup \partial A_2$}.
\eald \\
&\text{(solvability)}~~~ && 
\bald[t]
	&\text{If $\Deg(y;A;z)\neq 0$~ for a $z\notin y(\partial A)$}\\
	& \quad\text{then there exists $x\in A$ with $y(x)=z$.} 
\eald \\
&\bald[t] &\text{(homotopy}\\&~~\text{invariance)} \eald~~~&& \bald[t] 
&\Deg(y_t;A_t;z_t)=\Deg(y_0;A_0;z_0)~\text{for all $t\in [0,1]$,}\\
&\quad\text{if $z_t\notin y_t(\partial A_t)$ for all $t\in [0,1]$ and \eqref{homotopy} holds.}
\eald 
\ealdat
\]
Here, $y_t$ and $z_t$ are assumed to be a homotopies along $A_t$ in the sense that
\begin{align}\label{homotopy}
\bald
	&\text{$V:=\mysetr{(x,t)}{t\in [0,1]~\text{and}~x\in A_t}$}\\
	&\text{is bounded and open relative to $\RR^d\times [0,1]$,}\\
  &\text{$(t,x)\mapsto y_t(x),~~\bar{V}\to \RR^d$, is continuous and}\\
	&\text{$t \mapsto z_t,~~[0,1]\to \RR^d$, is continuous}.
\eald
\end{align}
In many cases, homotopy invariance is only stated and applied with cylinders $V=[0,1]\times A$, which also suffices for us here. For the general version see \cite{Kie12B} (e.g.).

Besides solvability and additivity, we here mainly use a few other properties of the degree which
can be derived from homotopy invariance:
\[
\baldat{2}
&\text{(continuity)}\quad && z\mapsto \deg(y;A;z)~~\text{is continuous on $\RR^d\setminus y(\partial A)$.} \\
&\text{(stability)}\quad && \bald[t]
	&\deg(y_1;A;z)=\deg(y_2;A;z) \\
	&\quad\text{if $\norm{y_1-y_2}_{C(\partial A;\RR^d)}<\dist{z}{y_1(\partial A)}$.} 
\eald\\
&\text{(bnd.~controlled)}\quad && \deg(y_1;A;\cdot)=\deg(y_1;A;\cdot)~~\text{if $y_1=y_2$ on $\partial A$.} \\
\ealdat
\]
Since the degree is integer-valued, continuity means it is locally constant. Stability is also a continuity property, now in $y$ instead of $z$.
"Boundary controlled" means that as far as $y$ is concerned, the degree is fully determined by the 
values of $y$ on $\partial A$. As we can always extend continuous functions from a compact set like $\partial A$ to the whole space, the degree is well-defined also for functions that are only given and continuous on $\partial A$.

To explicitly compute the degree in some examples, the following partial representation is helpful: If $y\in C^1(\bar A;\RR^d)$ and $z$ is a \emph{regular value} of $y$, i.e., 
$\det\nabla y(x)\neq 0$ for each $x\in y^{-1}(\{z\})$, then
\[
	\Deg(y;A;z)=\sum_{x\in y^{-1}(\{z\})} \operatorname{sgn}(\det\nabla y(x))\quad\text{if $z\notin y(\partial A)$}.
\]
Here, $\operatorname{sgn}$ denotes the sign ($\operatorname{sgn}(t)=t/\abs{t}$ if $t\neq 0$, and $\operatorname{sgn}(0)=0$).
This formula can also be used as the basis
of a definition of the degree. Besides, it 
determines the behavior of the degree with respect to reflections, a property which 
generalizes to all continuous $y$:
\[
	\Deg(Ry;A;Rz)=-\Deg(y;A;z) \quad\text{if $z\notin y(\partial A)$ and $R$ is a reflection},
\]
i.e., if $R\in O(d)$ with $\det R=-1$.

\section{The degree and approximate invertibility}\label{sec:AIBdeg}

Below, we repeatedly split sets in $\RR^d$ into their bounded and unbounded connected components. 
For this purpose, we introduce the following shorthand notation.
\begin{defn}[$\cB$ and $\cU$: bounded and unbounded components]\label{def:UB}~\\
Given a compact set $K\subset \RR^d$, we decompose 
\[
	\RR^d\setminus K~=~\cB(\RR^d\setminus K)~\cup~ \cU(\RR^d\setminus K),
\]
where $\cU(\RR^d\setminus K)$ denotes the unbounded connected component of $\RR^d\setminus K$, and $\cB(\RR^d\setminus K)$ denotes
the union of all bounded connected components of $\RR^d\setminus K$. If the choice of the set $K$ is obvious from the context,
we simply write $\cB$ and $\cU$.
\end{defn}

A core ingredient of this article is the following
statement about the degree of maps that are approximately invertible on the boundary.
\begin{thm}\label{thm:deg-of-def}
Let $\Omega\subset \RR^d$ 
be a bounded domain such that $\RR^d\setminus \partial\Omega$ has exactly two connected components, let
$y\in C(\bar\Omega;\RR^d)\cap \AIB$ and let
\[
	\imT(y;\Omega):=\mysetr{z\in \RR^d\setminus y(\partial\Omega)}{\Deg(y;\Omega;z)\neq 0}. 
\]
Then there exists a fixed $\sigma\in \{\pm 1\}$ 
such that
\begin{align}\nonumber
\hspace*{6ex}\bald
	&\Deg(y;\Omega;z)=\sigma \quad \text{for every $z\in \imT(y;\Omega)$}
		\eald
\end{align}
%\begin{align}\nonumber
%\hspace*{6ex}\bald
	%&\imT(y;\Omega):=\mysetr{z\in \RR^d\setminus y(\partial\Omega)}{\Deg(y;\Omega;z)\neq 0}.
	%%&\text{where }~\hat\cU:=\mysetr{z\in \RR^d\setminus y(\partial\Omega)}{\Deg(y;\Omega;z)=0}&.
	%\eald
%\end{align}
Moreover, $\imT(y;\Omega)$ is an open subset of $\cB(\RR^d\setminus y(\partial\Omega))$, and any connected component of $\cB(\RR^d\setminus y(\partial\Omega))$ 
is either fully contained in $\imT(y;\Omega)$ or does not intersect $\imT(y;\Omega)$ at all.
\end{thm}
\begin{rem}
Even if $y\in \AI(\bar\Omega)\subset \AIB$, it can happen that $\imT(y;\Omega)$ is not connected, for instance when $y$ compresses a surface bisecting $\Omega$ to a point.
\end{rem}
\begin{rem}
If $y|_{\partial\Omega}$ is already injective, it suffices to use $\varphi_k:=y|_{\partial\Omega}$ to show that $y\in\AIB$.
In this special case, Theorem~\ref{thm:deg-of-def} reduces to \cite[Proposition 2.2]{MueSpeTa96a}, 
which follows from the proof of a generalization of the Jordan Separation Theorem (see, e.g., \cite[Theorem 3.29 and its proof]{FoGa95B}). 
However, assuming invertibility on the boundary rules out deformations with
self-contact, while Theorem~\ref{thm:deg-of-def} still applies with a suitable sequence $\varphi_k$. 
\end{rem}
\begin{rem}
The assumption that $\RR^d\setminus \partial\Omega$ has only two connected components cannot be dropped, not even for orientation preserving maps. 
Counterexamples are given in Appendix~\ref{sec:A:topological}. Topologically speaking, this assumption 
expresses that $\partial\Omega$ is oriented in a degenerate sense inherited from $\Omega$ and the ambient space.
\end{rem}
\begin{rem}\label{rem:Lipbndcon}
Theorem~\ref{thm:deg-of-def} in particular applies for every bounded Lipschitz domain $\Omega$ with connected boundary. 
Here and throughout the article, Lipschitz domain is understood in the strong sense, i.e., the boundary can be locally (in $\RR^d$) represented as the graph of a Lipschitz function.
In this case, $\RR^d\setminus \partial\Omega$ has exactly two connected components, $\Omega$ being the bounded one, because 
a neighborhood of $\partial\Omega$ in $\Omega$ is connected.
\end{rem}

\begin{rem}[A degree on $y(\partial\Omega)$]
If $\imT(y;\Omega)\neq \emptyset$, Theorem~\ref{thm:deg-of-def} 
allows us to define a degree for continuous maps $f:y(\partial\Omega)\to y(\partial\Omega)$.
For that purpose, choose a point $z_0\in \imT(y;\Omega)$ (i.e., so that $\Deg(y;\Omega;z_0)=\sigma\in \{\pm 1\}$) and
a constant $\sigma_0\in \{\pm 1\}$
and define
\[
	\Deg_{y(\partial\Omega)}(f):=\sigma_0 \Deg(\hat{f}\circ y;\Omega;z_0), %=\sigma_0\lim_{k\to\infty}\Deg(\hat{f}\circ \hat{y}_k;\Omega;z_0),
\]
where $\hat{f}:\RR^d\to\RR^d$ is an arbitrary continuous extension of $f$.
This is well defined, since the right hand side only depends on the values of $\hat{f}\circ y|_{\partial\Omega}=f\circ y|_{\partial\Omega}$. 
The definition depends both on $y$ and on $z_0$. However, the dependence on $y$ (as a parametrization of $y(\partial\Omega)$) only enters through the orientation encoded in $\sigma$, and if we fix a connected component $\hat{B}$ of $\im_T(y;\Omega)$, the definition 
the definition does not depend to the choice of $z_0$ within $\hat{B}$.
In fact, 
we thus naturally obtain a whole family of degrees, one for each such connected component $\hat{B}$.
%Among other features, for invertible maps $f_1,f_2$, this degree inherits the multiplicative property 
%$\Deg_{y(\partial\Omega)}(f_1\circ f_2)=\Deg_{y(\partial\Omega)}(f_1)\Deg_{y(\partial\Omega)}(f_2)$
%from the corresponding one on $\Omega$ 
%which we exploit in the proof of Theorem~\ref{thm:deg-of-def}.
%As a consequence, we can always get that $\Deg_{y(\partial\Omega)}(\identity)=1$ by a choice of $\sigma_0$. 
 %All of the above is easy to see if $y|_{\partial\Omega}$ is invertible, 
%but it also holds in the more general case of $y\in \AIB$, where we can use the boundary approximations $y_k$.
\end{rem}

\begin{proof}[Proof of Theorem~\ref{thm:deg-of-def}]
Let
\[
	\hat\cU:=\mysetr{z\in \RR^d\setminus y(\partial\Omega)}{\Deg(y;\Omega;z)=0}=\RR^d\setminus [y(\partial\Omega)\cup \imT(y;\Omega)]
\]
(i) First observe that by the solvability property, $\Deg(y;\Omega;z)=0$ for some $z\in \cU:=\cU(\RR^d\setminus y(\partial(\Omega))$, since 
there clearly exist values $z\in \cU$ which are not in the bounded set $y(\Omega)$.
As $\Deg(y;\Omega;\cdot)$ is constant on connected components of $\RR^d\setminus y(\partial\Omega)$ (continuity), 
we infer that $\cU\subset \hat{\cU}$. In addition, 
connected components of $\cU$ are either fully contained in $\hat{\cU}$ or do not intersect $\hat{\cU}$ at all.

It remains to show that $\Deg(y;\Omega;z)=\sigma\in \{\pm 1\}$ for $z\in \imT(y;\Omega)$. By the Tietze extension theorem, both $\varphi_k$ and $\varphi_k^{-1}$ have continuous 
(but not necessarily invertible) extensions to $\RR^d$, say, $Y_k$ and $Z_k$. Clearly,
\[
	(Z_k\circ Y_k)(x)=x\quad\text{for every $x\in\partial\Omega$},
\]
and since the degree only depends on the values on the boundary, this implies that
\begin{align}\label{degcalc-1}
	\Deg(Z_k\circ Y_k;\Omega;x)=\Deg(\identity;\Omega;x)=\left\{\baldat[c]{2}
	&1\quad && \text{for every $x\in\Omega$,}\\
	&0\quad && \text{for every $x\in \RR^d\setminus\bar\Omega$.}
	\ealdat \right.
\end{align}
Due to our topological assumption on $\Omega$ and the Jordan Separation Theorem (see, e.g., \cite[Theorem 3.29 ]{FoGa95B}),
$\RR^d\setminus \varphi_k(\partial\Omega)$ also has exactly two connected components, 
$\cB_k:=\cB(\RR^d\setminus \varphi_k(\partial\Omega))$ and $\cU_k:=\cU(\RR^d\setminus \varphi_k(\partial\Omega))$.
Moreover, $\cB_k\subset Y_k(\Omega)$ and $\Deg(Y_k;\Omega;z)=0$ for $z\in \cU_k$.
As a consequence of the multiplication formula for the degree \cite[Theorem 2.10]{FoGa95B}, we get that
\begin{align}\label{degcalc-mult}
  \Deg(Z_k\circ Y_k;\Omega;x)=\Deg(Z_k;\cB_k;x) \Deg(Y_k;\Omega;z)\quad\text{for $z\in\cB_k$.}
\end{align}
Here, notice that $\Deg(Y_k;\Omega;\cdot)$ is constant on $\cB_k$. In view of \eqref{degcalc-1} and the fact that the degree is integer-valued,
\eqref{degcalc-mult} entails that for each $k$, there exists a $\sigma_k\in\{\pm 1\}$ such that
\begin{align}\label{degcalc-2} 
  \sigma_k=\Deg(Z_k;\cB_k;x)=\Deg(Y_k;\Omega;z)\quad\text{for $x\in\Omega$ and $z\in\cB_k$.}
\end{align}
In addition, by the stability of the degree,
\begin{align}\label{degcalc-3}
\bald
  &\Deg(Y_k;\Omega;z)= \Deg(y;\Omega;z)\quad \\
	&\text{for every $z\in \RR^d$ with }\dist{z}{y(\partial\Omega)}>\norm{Y_k-y}_{C(\partial\Omega)}.
\eald
\end{align}
Since $\Deg(y;\Omega;z)=0$ for all $z\in\hat{\cU}$,
\eqref{degcalc-2} and \eqref{degcalc-3} imply that
\begin{align}\label{degcalc-4a}
  \cB_k\supset V_k:=\imT(y;\Omega)\cap\mysetr{z}{
	\begin{array}{l}
		\text{$z\in \cB(\RR^d\setminus y(\partial\Omega))$ and}\\
		\dist{z}{y(\partial\Omega)}>\norm{Y_k-y}_{C(\partial\Omega)}
	\end{array}
	}
\end{align}
and
\begin{align}\label{degcalc-4}
	\sigma_k=\Deg(y;\Omega;\cdot) \quad\text{is constant on $V_k$ for each $k$}.
\end{align}
Since $\norm{Y_k-y}_{C(\partial\Omega)}\to 0$ as $k\to\infty$,
the limit of the increasing sequence of open sets $V_k$ is given by
\begin{align}\nonumber%\label{degcalc-5}
\bald
	\imT(y;\Omega)=\bigcup_k V_k & = \mysetr{z\in\RR^d}{
	\begin{array}{l}
		\text{$z$ is in a bounded connected}\\
		\text{component of $\RR^d\setminus (y(\partial\Omega)\cup \cU)$}
	\end{array}
	}\\
	& \subset \mysetr{z\in \RR^d\setminus y(\partial\Omega)}{z\in \cB_k~~\text{for 
	almost all $k$}}.
\eald
\end{align}
By \eqref{degcalc-2} and \eqref{degcalc-3}, we also get that $\sigma_k\to \sigma$, whence
\begin{align}\label{degcalc-6}
\bald
  &\Deg(y;\Omega;z)=\sigma \quad 
	\text{for all $z\in \imT(y;\Omega)$, with fixed $\sigma\in\{\pm 1\}$}.
\eald
\end{align}
The remaining assertions are a straightforward consequence of the fact that $\Deg(y;\Omega;\cdot)$ is constant on each connected component of 
$\RR^d\setminus y(\partial\Omega)$.
\end{proof}

The equivalent of Theorem~\ref{thm:deg-of-def} for $\AI(\bar\Omega)$ is a similar but easier because the degree of full homeomorphisms is known. It works for all domains:
\begin{thm}\label{thm:deg-of-def2}
Let $\Omega\subset \RR^d$ 
be a bounded domain and let $y\in C(\bar\Omega;\RR^d)\cap \AI(\bar\Omega)$.
Then all the assertions of Theorem~\ref{thm:deg-of-def} hold.
%Moreover, $\imT(y;\Omega)$ is an open subset of $\cB(\RR^d\setminus y(\partial\Omega))$, and any connected component of $\cB(\RR^d\setminus y(\partial\Omega))$ 
%is either fully contained in $\imT(y;\Omega)$ or does not intersect $\imT(y;\Omega)$ at all.
\end{thm}
\begin{proof}
For the approximating homeomorphisms $\varphi_k:\bar\Omega\to \varphi_k(\bar\Omega)$,
it is well known that $\Deg(\varphi_k;\Omega;z)=\sigma_k\in \{-1,1\}$ for every $z\notin \varphi_k(\partial\Omega)$
(this is a simpler application of the multiplication theorem for the degree), 
and that the bounded connected component of $\RR^d\setminus \varphi_k(\partial\Omega)$ is given by $\cB_k:=\varphi_k(\Omega)$.
The rest of the proof is analogous to the proof of Theorem~\ref{thm:deg-of-def}.
%Since $\varphi_k\to y$ uniformly on $\bar\Omega$, the stability of the degree implies that
%$\sigma_k$ must converge to some limit $\sigma\in \{-1,1\}$,
%and $\sigma=\Deg(y;\Omega;z)$ for all $z\in \imT(y;\Omega)$.
%Here, notice that for $z\in \imT(y;\Omega)$, we have $\dist{z}{y(\partial\Omega)}>0$, and therefore
%$z\notin \varphi_k(\partial\Omega)$ for $k$ large enough.
%
%As before, the remaining assertions are a straightforward consequences of properties of the degree.
\end{proof}

\section{The degree and orientation preserving maps}\label{sec:top}

In this section, we study the consequences of Theorem~\ref{thm:deg-of-def} for the case of strictly orientation preserving maps.
Some of the results in this section are essentially known (cf.~\cite{Sve88a}, in particular) and the technique is classical, but there is no comprehensive collection suitable for our purposes.
For this reason, full proofs are given throughout. To the best of my knowledge, the global results summarized Theorem~\ref{thm:sop-deg1} are new in this generality, in particular the description of sets mapped to a point stated in terms that work well in $\DEGoneloc$ as illustrated by Corollary~\ref{cor:sop-deg1loc}. Already in \cite{Ba81a}, essentially the same properties were obtained, but only for deformations whose boundary values are continuous and admit a homeomorphic extension, which among other things provides a much more straightforward way to control of $y(\bar\Omega)$, whereas we cannot even assume that $y(\partial\Omega)$ is defined. 

Proofs of partially related results in Nonlinear Elasticity, notably \cite{Ba81a}, \cite{Sve88a}, \cite{MueSpe95a} and \cite{MueSpeTa96a}, the degree-theoretic parts of \cite{HeKo14B} and, more recently, \cite{BaHeMo17a}, mix topological arguments with fine properties of Sobolev functions. Here, the presentation is purely topological, including a notion of strictly orientation preserving maps based on the degree. It expands and complements some results of \cite{TiYou62a},
and it also reveals that some of the auxiliary results of \cite{Sve88a} and related papers can in fact be derived using only information about the degree. 

We start with general properties of orientation preserving maps, defined topologically in terms of the degree. As illustrated by \cite{BaHeMo17a}, this is also useful in settings with low Sobolev regularity, suitable subsets of $W_+^{1,p}$ with $d-1<p<d$. 
Such settings as treated in \cite{Sve88a}, \cite{MueSpe95a}, \cite{HeMo10a} and \cite{BaHeMo17a} still allow degree theory, but only for selective "good" open subsets of $\Omega$. In our notation below, $\cO$ could actually play that role. 
However, the more difficult issue in this scenario is that $y$ than can only be assumed to be continuous on boundaries of good sets, not everywhere. For global theory, using only balls as good sets (like in \cite{Sve88a} or \cite{MueSpe95a}) does not suffice, because they in general cannot separate connected components. In any case, in this article, we will not pursue this topic further.

\subsection{Strictly orientation preserving maps}  

\begin{defn}[non-degenerate, (strictly) orientation preserving]\label{def:op}~\\
Let $\Omega\subset \RR^d$ be open, let $y:\Omega\to\RR^d$ be continuous and let 
\[
	\cO:=\{A\subset \Omega\mid \text{$A$ is open and $A\subset \subset\Omega$}\}.
\]
We say that $y$ is \emph{non-degenerate} if for each non-empty $U\in \cO$, 
there exist $A\in \cO$ with $A\subset U$ and $z\in \RR^d\setminus y(\partial A)$ with $\Deg(y;A;z)\neq 0$.
We say that $y$ is \emph{orientation preserving}
if $\Deg(y;\cdot;\cdot)\geq 0$, i.e., $\Deg(y;A;z)\geq 0$ for all $A\in \cO$ and all $z\in\RR^d\setminus y(\partial A)$.
If $y$ is non-degenerate and orientation preserving
%and for every non-empty open $A\subset\subset \Omega$, 
%there exists $x\in A$ such that $y(x)\notin y(\partial A)$ and $\Deg(y;A;y(x))>0$,
then $y$ is called \emph{strictly orientation preserving}.
%In this case, we formally write $\Deg(y;\cdot;\cdot)>0$. 
\end{defn}
For instance, any homeomorphism $y$ 
is non-degenerate and has constant degree in $\{\pm 1\}$ on its image.
Up to a reflection to achieve degree $+1$ if necessary, it is also strictly orientation preserving.
\begin{rem}
The same topological notion of orientation preserving maps
is used in \cite{BaHeMo17a}. By itself, it is often too weak to be useful because it allows examples where the degree
is simply not defined on the image of $y$ on all of $\Omega$ or large subsets.
The attributes "non-degenerate" and "strictly orientation preserving" as defined here are not standard. The latter is similar but not equivalent to "locally sense-preserving" in the sense of, e.g., \cite{TiYou62a} (cf.~Remark~\ref{rem:TiYo}). "Non-degenerate" essentially expresses that 
$y$ is not allowed to compress open sets to sets with empty interior. 
%Here, notice that by the solvability property of the degree, $\Deg(y;A;z)\neq 0$ implies that 
%$y^{-1}(\{z\})\cap A\neq \emptyset$.
\end{rem}
\begin{rem}
As we will see later, for $p\geq d$, deformations $y\in W_{+}^{1,p}$ are always strictly orientation preserving in the sense of Definition~\ref{def:op}. The case $p=d$ 
is the main reason we want $\cO$ to be the family of open sets compactly contained in $\Omega$ instead of all open $A\subset \Omega$. As a consequence, it is not needed to have $y$ continuous up to the boundary. 
\end{rem}

For strictly orientation preserving maps in the sense of Definition~\ref{def:op}, both the local  and the global degree are actually positive for any admissible value in the image of $y$:
\begin{lem}\label{lem:sop}
Let $\Omega\subset\RR^d$ be open, let $\cO=\mysetr{A\subset\subset\Omega}{\text{$A$ is open}}$, let $y\in C(\Omega;\RR^d)$ be strictly orientation preserving in the sense of Definition~\ref{def:op}, and
let $z\in y(\Omega)$.  Then for all connected components $C_z$ of $y^{-1}(\{z\})$ with $C_z\subset\subset \Omega$, there exists a sequence 
of sets $A_n \in \cO$ such that for every $n\in \NN$, $A_{n+1}\subset A_n$,
\begin{align}\label{lemop-31}
\bald
	&C_z\subset A_n,~~ \Dist{C_z}{\partial A_n}\leq \tfrac{1}{n},\\
	&z\notin y(\partial A_n)~~\text{and}~~\Deg(y;A_n;z)\geq 1.
\eald
\end{align}
If we also have that $\Omega$ is bounded and $y\in C(\bar\Omega;\RR^d)$, then 
\[
	\text{$\Deg(y;\Omega;z)\geq 1$ for all $z\in y(\Omega)\setminus y(\partial\Omega)$}.
\]
\end{lem}
\begin{proof}
Since $C_z$ is a compact subset of $\Omega$,
it has a positive distance to $\partial\Omega$. 
Choose a sequence of open sets $U_n\in\cO$ with 
\begin{align}\label{lemop-30a}
	C_z\subset U_n\subset\subset\Omega\quad\text{and}\quad\Dist{C_z}{\partial U_n}\leq \tfrac{1}{n}.
\end{align}
The set 
\[
	S_n:=\bar{U}_n\cap y^{-1}(\{z\})
\]
is compact. If 
$S_n=C_z$ then $y^{-1}(\{z\})\cap \partial U_n=\emptyset$ and we take $A_n:=U_n$. Otherwise,
$S_n$ is not connected
and can be 
separated into two disjoint compact subsets $T_n$ and $R_n=S_n\setminus T_n$ such that
\begin{align}\label{lemop-30b}
	C_z\subset T_n\quad\text{and}\quad S_n\cap \partial U_n\subset R_n
\end{align}
(Whyburn Lemma, see \cite[(9.3) in Chap.~I, p.~12]{Why58B}).
Since $T_n$ and $R_n$ have positive distance,
there exists $A_n\in \cO$, $A_n\subset U_n$ with $C_z\subset T_n\subset A_n$ and  $R_n\cap \partial A_n=\emptyset$.
In all cases, we conclude that
\begin{align}\label{lemop-32}
	C_z\subset A_n,~~~ \Dist{C_z}{\partial A_n}\leq \frac{1}{n}~~~\text{and}~~~z\notin y(\partial A_n).
\end{align}
W.l.o.g., we may also assume that $A_{n+1}\subset A_n$ for all $n$.
Since $y$ is strictly orientation preserving, there exists
\begin{align}\label{lemop-33}
	\text{$x_n\in A_n$ with $y(x_n)\notin \partial A_n$ and $\Deg(y;A_n;y(x_n))\geq 1$,}
\end{align}
the latter because $\Deg(y;A_n;y(x_n))>0$ and the degree is integer-valued.
Due to \eqref{lemop-32}, we also have that $\dist{x_n}{C_z}\to 0$, and thus $y(x_n)\to z$ by the (locally uniform) continuity of $y$.
Using the additivity and continuity of the degree together with $\Deg(y;\cdot;\cdot)\geq 0$ ($y$ is orientation preserving), \eqref{lemop-33}
implies that for all $n$ and every $k\geq n$ large enough so that $y(x_k)\notin \partial A_n$ 
(which is possible since $y(x_k)\to z$ and $z$ has a positive distance to the compact set $\partial A_n$ for fixed $n$),
\[
\bald
	1 \leq \Deg(y;A_k;y(x_k)) \leq \Deg(y;A_n\setminus \bar{A}_k;y(x_k))+ \Deg(y;A_k;y(x_k)) &\\
	 =\Deg(y;A_n;y(x_k)) \underset{k\to\infty}{\To} \Deg(y;A_n;z)&. 
\eald
\]
We conclude that $\Deg(y;A_n;z)\geq 1$ for all $n$. Combined with \eqref{lemop-32}, 
this concludes the proof of \eqref{lemop-31}.

Finally, if $y$ is continuous up to the boundary and  $z\in y(\Omega)\setminus y(\partial\Omega)$, then $y^{-1}(\{z\})\cap \Omega\neq \emptyset$ and
any of its connected components $C_z$ is compactly contained in $\Omega$. Thus, $\Deg(y;\Omega;z)\geq \Deg(y;A_n;z)\geq 1$ by additivity of the degree.
\end{proof}

\begin{lem}[on topological image and reduced domain]\label{lem:orpres-top} ~\\
Let $\Omega\subset\RR^d$ be open, let $\cO=\mysetr{A\subset\subset\Omega}{\text{$A$ is open}}$, let $y\in C(\Omega;\RR^d)$, 
let $A,A_1,A_2\in \cO$ and let $U\subset \Omega$ be open.
We define the open sets
\begin{alignat*}{2}
	&\imT(y;A)&&:=\mysetr{z\in\RR^d\setminus y(\partial A)}{\Deg(y;A;z)\neq 0}, \\
	&\imloc(y;U)&&:= \textstyle{\bigcup_{A\in \cO,A\subset U}} \imT(y;A) 
\end{alignat*}
and the associated \emph{(topologically) reduced domain}
\begin{alignat*}{2}
	&\cR_y(U)&&:=\mysetb{x\in U}{\exists A\in \cO:~x\in A\subset U~\text{and}~y(x)\notin y(\partial A)}.
\end{alignat*}
Then we have the following:
\begin{itemize}
\item[(i)] $\imT(y;A)\subset \interior y(A)$, the interior of $y(A)$.
%\hat\cB\subset y(\Omega)\setminus y(\partial\Omega)$.
\item[(ii)] If $y$ is orientation preserving then 
\[
		A_1\subset A_2 ~~~\Longrightarrow~~~ \imT(y;A_1) \subset \imT(y;A_2)\cup y(\partial A_2).
\]
%$\Deg(y;\Omega;\cdot)\geq 0$.
\item[(iii)]
If $y$ is strictly orientation preserving, then 
\[
\bald
	&\imT(y;A)=y(A)\setminus y(\partial A), \\
	&\imloc(y;U)\subset y(U) \subset \overline{\imloc(y;U)}, \\
	&\text{$U\setminus \cR_y(U)$ has empty interior,} \\
	&y(\cR_y(U))=\imloc(y;U),~\text{and}\\
	&y(\bar{U}\setminus \cR_y(U))\subset y(\partial U).
\eald
\]
\end{itemize}
If we have in addition that $\Omega$ is bounded and $y\in C(\bar\Omega;\RR^d)$, then all of the above also holds
for $A=A_2=U=\Omega$.
\end{lem}
\begin{rem}\label{rem:imT}
In context of Nonlinear Elasticity, the concept of the topological image $\imT$ (without the name) was 
introduced in \cite{Sve88a}. We use it here without artificially adding the image of the boundary,
which makes it an open set (due to continuity of the degree). This has the disadvantage that in general,
full monotonicity with respect to the domain cannot be expected for even in case of strictly orientation preserving maps, only what follows from Lemma~\ref{lem:orpres-top} (ii) and (iii).
The "localized" topological image $\imloc(y;U)$ fixes this issue by collecting all points in local topological images. A variant of it was also used in 
\cite{Ta88a}. 

As shown in Theorem~\ref{thm:sop-deg1} below, for strictly orientation preserving maps with global degree $\leq 1$, $\imloc(y;U)$ coincides with $\imT(y;U)$ as long as $y(\partial U)$ has empty interior.
In addition, $\imloc(y;U)$ is a natural substitute of $y(U)\setminus y(\partial U)$.
Unlike the latter, it remains meaningful with $U=\Omega$ even if $y$ cannot be continuously extended to $\partial\Omega$,
and it possibly contains more information in cases where $y|_{\partial U}$ exhibits wild behavior like Peano curves. 
\end{rem}
\begin{rem}\label{rem:Ry}
As Lemma~\ref{lem:orpres-top} and Theorem~\ref{thm:sop-deg1} show, the topologically reduced domain $\cR_y$ 
introduced in the lemma is always naturally associated to $\imloc(y;\cdot)$. By its definition and the continuity of $y$, 
$\cR_y(U)$ is open. Moreover, since 
we can always isolate connected components with an arbitrarily close surrounding neighborhood unless they touch the boundary (cf.~Lemma~\ref{lem:sop}), 
\begin{align}\label{Ry-repres1}
	\cR_y(U)=\Lambda:=\bigcup_{z\in \RR^d,~ C_z\in \operatorname{ICC}(z,U)} C_z,
\end{align}
where $\operatorname{ICC}(z,U)$ denotes the family of all "inner" connected components of $y^{-1}(\{z\})\cap U$,
i.e.,
\[
	\operatorname{ICC}(z,U):=\mysetr{C_z}{
	\begin{array}{l}
	\text{$C_z$ is a connected component of $y^{-1}(\{z\})\cap U$}\\
	\text{with $C_z\subset\subset U$}
	\end{array}
	}
\]
This also entails that for every $z\in y(\partial U)$ and every component $C_z$ of $y^{-1}(\{z\})\cap U$ with $C_z\not\subset \cR_y(U)$,
we have that $C_z\cap \cR_y(U)=\emptyset$ and $\bar{C}_z\cap \partial U\neq \emptyset$.
Interestingly, it is not obvious that $\Lambda$ is open as defined.
\end{rem}
\begin{rem}\label{rem:Ry2}
Another basic property of $\cR_y$ to keep in mind when working in $\DEGoneloc$ is that 
for every covering $(\Omega_m)_{m\in \NN}$ of $\Omega$ with open sets $\Omega_m\subset \Omega_{m+1}\subset\subset \Omega$,
\begin{align}\label{Ry-repres2}
\cR_y(\Omega_m)\subset \cR_y(\Omega_{m+1})\quad\text{and}\quad\cR_y(\Omega)=\textstyle \bigcup_{m\in \NN} \cR_y(\Omega_m),
\end{align}
a straightforward consequence of the definition of $\cR_y$. Similarly,
\begin{align}\label{imloc-repres2}
	\imloc(y;\Omega_m)\subset \imloc(y;\Omega_{m+1})\quad\text{and}\quad\imloc(y;\Omega)=\textstyle \bigcup_{m\in \NN} \imloc(y;\Omega_m).
\end{align}
\end{rem}
\begin{rem}[Connections to the results of Titus and Young \cite{TiYou62a}]\label{rem:TiYo}
The set $\cR_y(\Omega)$ is related to concepts appearing in \cite{TiYou62a}.
If we knew that
\begin{align}\label{rTY-1}
	\text{$y(\Omega\setminus \cR_y(\Omega))$ is closed with empty interior}
\end{align}
for a strictly orientation preserving $y$, then
$y$ would belong to the class of functions called $\Omega$ in \cite{TiYou62a} 
with the corresponding set $C_f:=\Omega\setminus \cR_y(\Omega)$ for $f=y$ (and our domain $\Omega$).
%(which originally was probably intended to be the set where $\det\nabla y=0$ in smoother scenarios). 
However, unless $\cR_y(\Omega)=\Omega$, \eqref{rTY-1} is extremely unnatural even if $y$ is continuous up to the boundary and $y(\partial\Omega)$ has empty interior; as a matter of fact,
it even implies that $\cR_y(\Omega)=\Omega$ in our setting.
In essence, the issue reflects that \cite{TiYou62a} was written for manifolds without boundary.
Nevertheless, the results of \cite{TiYou62a} do apply to $y|_{\cR_y(\Omega)}$ 
(with $f=y$, $N=\cR_y(\Omega)$, $M=\RR^d$ and $C_f=\emptyset$).
As is, \cite{TiYou62a} effectively cannot provide any information on $\Omega\setminus \cR_y(\Omega)$ without additional assumptions.
On such assumption is that $y$ is \emph{light}, i.e, $y^{-1}(\{z\})$ is totally disconnected for all $z$, but this again implies that $\Omega=\cR_y(\Omega)$.

By contrast, Theorem~\ref{thm:sop-deg1} 
and Corollary~\ref{cor:sop-deg1loc}
discuss
$y$ on $\Omega\setminus \cR_y(\Omega)$, and we already know that $\Omega\setminus \cR_y(\Omega)$ has empty interior by Lemma~\ref{lem:orpres-top} (iii) (and Remark~\ref{rem:Ry2}), essentially a consequence of our stronger notion of strictly orientation preserving maps.
Unlike the latter, the notion of sense preserving used in \cite{TiYou62a} never holds at points $x\in\Omega\setminus \cR_y(\Omega)$ because it requires the connected component $C_x$ of $x$ in $y^{-1}(\{y(x)\})$ to be compact in $\Omega$ (in other words, $C_x$ must not touch $\partial\Omega$). On $\cR_y(\Omega)$, it is weaker than "strictly orientation preserving", a consequence of Lemma~\ref{lem:sop}.
\end{rem}

\begin{proof}[Proof of Lemma~\ref{lem:orpres-top}]
Notice that $\imT(y;A)$ is open by the continuity of the degree. Consequently, $\imloc(y;U)$ is open as the union of open sets.

{\bf (i)} This is clear since $\imT(y;A)$ is open and $\imT(y;A)\subset y(A)$ by the solvability property of the degree.

{\bf (ii)} Let $z\in \imT(y;A_1)$. We may assume that $z\notin y(\partial A_2)$, because otherwise the assertion is obvious.
By definition of $\imT(y;A_1)$, $z\notin y(\partial A_1)$ and $\Deg(y;A_1;z)\neq 0$. Since $y$ is orientation preserving, 
we even have that $\Deg(y;A_1;z)>0$, an we also know that $\Deg(y;A_2\setminus \bar{A}_1;z)\geq 0$. By additivity of the degree,
this implies that
\[	
 \Deg(y;A_2;z)=\Deg(y;A_2\setminus \bar{A}_1;z_0)+\Deg(y;A_1;z)>0.
\]
Consequently, $z\in \imT(y;A_2)$.

{\bf (iii)} {\bf "$\mbf{\imT(y;A)=y(A)\setminus y(\partial A)}$"}: "$\subset$" is clear by definition, and "$\supset$" is a consequence of Lemma~\ref{lem:sop}.

{\bf "$\mbf{\imloc(y;U)\subset y(U)}$"}: This follows from (i).

{\bf "$\mbf{y(U) \subset \overline{\imloc(y;U)}}$"}:  Let $z\in  y(U)$, and choose a point $x \in U$ with $y(x)=z$. Since $U$ is open, there exists positive radii $r(k)\to 0$ such that $B_{r(k)}(x)\subset\subset U$ for each $k$. As $y$ is strictly orientation preserving, there exist open $A_k\subset B_{r(k)}(x)$ and $x_k\in A_k$ with $y(x_k)\notin y(\partial A_k)$ and $\Deg(z;A_k;y(x_k))\neq 0$. By construction of the sets $A_k$, we have that $x_k\to x$ as $k\to\infty$. Consequently, $z=y(x)=\lim_k y(x_k)\in \overline{\imloc(y;U)}$.

"{\bf $\mbf{U\setminus \cR_y(U)}$ has empty interior}": Let $V\subset U\setminus \cR_y(U)$ be open. 
By definition of $\cR_y(U)$, we infer that
$y(x)\in y(\partial A)$ for all $x\in V$ and all $A\in \cO$ with $x\in A\subset U$. But since $y$ is strictly orientation preserving, only empty $V$ can satisfy this.

"$\mbf{y(\cR_y(U))=\imloc(y;U)}$": For every $x\in \cR(y;U)$ choose an admissible associated $A\in \cO$ from its definition,
and let $C_x$ denote the connected component of $y^{-1}(y(x))$ containing $x$. Since 
$y(x)\notin y(\partial A)$, $C_x$ is a compact subset of $A$. By Lemma~\ref{lem:sop}, there exists 
$A_1\in \cO$ with $C_x\subset A_1\subset A$, $y(x)\notin y(\partial A_1)$ and $\Deg(y;A_1;y(x))\geq 1$. Hence,
$y(x)\in \imT(y;A_1)\subset \imloc(y;U)$. Conversely, if $z\in \imloc(y;U)$, then $z\in \imT(y;A)$ for an $A\in \cO$, $A\subset U$ with $z\notin y(\partial A)$. In addition, $\Deg(y;A;z)\neq 0$, and consequently, there exists $x\in A$ with $y(x)=z$. We infer that $x\in \cR(y;U)$ and therefore $z=y(x)\in y(\cR(y;U))$.

"$\mbf{y(\bar U\setminus \cR_y(U))\subset y(\partial U)}$": It suffices to show that $\bar{U}\setminus y^{-1}(y(\partial U))\subset \bar{U}\setminus \cR_y(U)$.
Consider an arbitrary $x\in \bar{U}$ with $y(x)\notin y(\partial U)$. Then the connected component $C_x$ of $x$ in $\bar{U}\cap y^{-1}(y(x))$ has positive distance to $\partial U$. By Lemma~\ref{lem:sop}, there exists $A\in \cO$ with $C_x\subset A\subset U$ and $y(x)\notin y(\partial A)$. Hence, $x\in \cR_y(U)$, and since $\cR_y(U)$ is open, we even get that $x\notin \bar{U}\setminus \cR_y(U)$.  

{\bf  The case $\mbf{y\in C(\bar\Omega;\RR^d)}$}: In this case, the proofs of (i) and (iii) still work for $A=U=\Omega$. 
The same is true concerning (ii), as long as we can still show that $\Deg(y;\Omega\setminus \bar{A}_1;z)\geq 0$ for $z\notin y(\partial\Omega)$.
Here, notice that $A:=\Omega\setminus \bar{A}_1\notin \cO$ is not admissible in our definition of orientation preserving maps.
Instead, choose a sequence $\Omega_m\in \cO$ such that $\bar\Omega_m\subset \Omega$ and $\Dist{\Omega_m}{\partial\Omega}\to 0$ as $m\to\infty$.
Since $y$ is uniformly continuous, $z\notin y(\bar\Omega \setminus \Omega_m)$ for $m$ large enough. 
By additivity of the degree, this implies that $\Deg(y;\Omega\setminus \bar{A}_1;z)=\Deg(y;\Omega_m \setminus \bar{A}_1;z)\geq 0$. 
The inequality holds because $y$ is orientation preserving and $\Omega_m \setminus \bar{A}_1\in \cO$.
\end{proof}

\subsection{Strictly orientation preserving maps of degree one}  

The following theorem is of particular interest from the point of view of Nonlinear Elasticity, 
because it provides a major step towards the invertibility of deformations.
It summarizes the topological properties of strictly orientation preserving maps with degree $\leq 1$.
Among other things, it asserts that they are essentially monotone 
in the topological sense (see (ii) below). The only possible case where the preimage of a value can have
several connected components occurs when all of these are contracted by the deformation to boundary points where the deformed configuration is in self-contact (cf.~(iii)). 
These are exactly the pieces missing in the (topologically) reduced domain $\cR_y$ defined and studied in Lemma~\ref{lem:orpres-top}, cf.~Remark~\ref{rem:Ry}. 
\begin{thm}\label{thm:sop-deg1}
Let $\Omega\subset\RR^d$ be open, let $\cO=\mysetr{A\subset\subset\Omega}{\text{$A$ is open}}$, let $y\in C(\Omega;\RR^d)$ be strictly orientation preserving in the sense of Definition~\ref{def:op},
let $U\in \cO$ and let $\imT$, $\imloc$ and $\cR_y$ be defined as in Lemma~\ref{lem:orpres-top}.
In addition, assume that $\Deg(y;U;\cdot)\leq 1$. Then we have the following:
\begin{itemize}
\item[(i)] Let $A\in \cO$ with $A\subset U$. Then $\Deg(y;A;\cdot)\geq 1$ on $y(A)\setminus y(\partial A)$,
and $\Deg(y;A;\cdot)\leq 1$ on $y(A)\setminus [y(\partial A)\cup y(\partial U)]$. 
%IF POSSIBLE AT ALL, THIS NEEDS ARGUMENTS OF (ii) AND (iii):
%If $A\subset \cR_y(U)$ and $\partial A$ has empty interior, then $y(\partial A)$ has empty interior and
%$y(A)\setminus y(\partial A)\neq \emptyset$.
If $y(\partial U)$ has empty interior, then
$\Deg(y;A;\cdot)=1$ on $y(A)\setminus y(\partial A)$.
\item[(ii)] For all $z\in \imloc(y;U)\setminus y(\partial U)$,
\[
	\cR_y(U) \cap y^{-1}(\{z\})~~\text{is connected.}
\]
\item[(iii)] If $y(\partial U)$ has empty interior, then 
\[
\bald
	&\imloc(y;U)\cap  y(\partial U)=\emptyset. 
\eald
\]
The latter implies that
\[
\bald
&\text{$\cR_y(U)=U\setminus y^{-1}(y(\partial U))$,} \\
&\text{$\imloc(y;U)=y(U)\setminus y(\partial U)=\imT(y;U)$ and}\\ 
&\text{$y(\partial U)=\partial (\imloc(y;U))$.}
\eald
\]
\end{itemize}
If we have in addition that $\Omega$ is bounded and $y\in C(\bar\Omega;\RR^d)$, then we may also use $U=\Omega$.
\end{thm}
If we can find one $U$ such that $y(\partial U)$ has empty interior, then smaller sets with a reasonable boundary inherit this property -- a topological analogue of Lusin's property (N):
\begin{cor}\label{cor:sop-deg1}
Under the assumptions of Theorem~\ref{thm:sop-deg1}, 
suppose in addition that $\interior y(\partial U)=\emptyset$. Then 
\begin{align}\label{topLusin}
	\interior K=\emptyset\quad \Longrightarrow \quad \interior y(K)=\emptyset,\quad
	\text{for all compact $K\subset \bar{U}$.}
\end{align}
In particular, for any open and bounded $\tilde{U}\subset U$ with $\interior \partial \tilde{U}=\emptyset$,
all assertions in Theorem~\ref{thm:sop-deg1} (i) and (iii) hold for $\tilde{U}$ instead of $U$.

If $\Omega$ is bounded and $y\in C(\bar\Omega;\RR^d)$, then $U=\Omega$ is also admissible. 
\end{cor}
The next corollary summarizes our results for the case $y\in \DEGoneloc$. We now also assume that $\Omega$ is bounded to avoid complications that would appear otherwise.
\begin{cor}\label{cor:sop-deg1loc}
Let $\Omega\subset\RR^d$ be open and bounded, let $y\in C(\Omega;\RR^d)$ be strictly orientation preserving in the sense of Definition~\ref{def:op}, let $\cO=\mysetr{A\subset\subset\Omega}{\text{$A$ is open}}$
and let $\imT$, $\imloc$ and $\cR_y$ be defined as in Lemma~\ref{lem:orpres-top}.
In addition, assume that $y\in \DEGoneloc$ with respect to a regular inner covering $(\Omega_m)_{m\in\NN}$ of $\Omega$, and that $y(\partial\Omega_m)$ has empty interior for all $m$. 
Then the following holds:
\begin{itemize}
\item[(i)] For all $A\in \cO$, $\Deg(y;A;\cdot)=1$ on $y(A)\setminus y(\partial A)$.
%If $\partial A$ has empty interior, then so has $y(\partial A)$, and
%$y(A)\setminus y(\partial A)\neq \emptyset$.
%
\item[(ii)] For all $z\in y(\cR_y(\Omega))$, 
$\cR_y(\Omega) \cap y^{-1}(\{z\})~~\text{is connected}$ and 
compactly contained in $\cR_y(\Omega)$.
\item[(iii)] For all $z\in y(\Omega)\setminus y(\cR_y(\Omega))$
and every connected component $C_z$ of $y^{-1}(\{z\})$,
$\bar{C}_z\cap \partial\Omega\neq \emptyset$.
\item[(iv)] $\imloc(y;\Omega)=y(\cR_y(\Omega))$ and $\Omega\setminus \cR_y(\Omega)$ has empty interior.
\end{itemize}
If we have in addition that $y\in C(\bar\Omega;\RR^d)$, then 
$\imloc(y;\Omega)=y(\Omega)\setminus y(\partial\Omega)$ and 
$\deg(y;\Omega;\cdot)=1$ on $\imloc(y;\Omega)$. 
\end{cor}
\begin{rem}[Interface to the analytical setting] \label{rem:interface} 
%Theorem~\ref{thm:sop-deg1} 
Corollary~\ref{cor:sop-deg1loc} is meant to be applied to functions $y\in \cY$ in the effective admissible set $\cY\subset \{E<\infty\}$ of a variational minimization problem 
for a functional $E$, for instance a nonlinear elastic energy set up in a Sobolev space
(e.g., $\cY\subset W^{1,p}_+$ for $p\geq d$).
%Typically, we choose a regular inner covering $(\Omega_m)_{m\in\NN}$ of $\Omega$ (cf.~Definition~\ref{def:regcover}), apply the theorem with $U:=\Omega_m$ for each $m$ and exploit Remark~\ref{rem:Ry2}. If $\Omega$ is bounded and $y\in C(\bar\Omega;\RR^d)$, $U=\Omega$ suffices instead.
The following ingredients are needed to make this work:
\begin{enumerate}
\item[(a)] An embedding of $\cY$ into $C(\Omega;\RR^d)$ (cf.~Remark~\ref{rem:pequald});
\item[(b)] a global invertibility constraint or a boundary condition which implies $\cY\subset \DEGoneloc$ (cf.~Remark~\ref{rem:invcompare});
\item[(c)] $y$ is strictly orientation preserving for $y\in \cY$
(cf.~Lemma~\ref{lem:orpres-W1p});
\item[(d)] a way to prove that 
$y(\partial \Omega_m)$ has empty interior for each $m$.
Since $\cL^d(\partial\Omega_m)=0$ by our definition of a regular inner covering,
Lusin's condition (N) is sufficient (cf.~Remark~\ref{rem:pequald}). 
\end{enumerate}
If $y\in C(\bar\Omega;\RR^d)$, 
then for (d), by Corollary~\ref{cor:sop-deg1}, it also suffices to have that $y(\partial\Omega)$ has empty interior.
\end{rem}
\begin{rem} \label{rem:globalneeded}
Without the global invertibility constraint $\deg(y;U;\cdot)\leq 1$, the assertions of Theorem~\ref{thm:sop-deg1} 
can fail to hold even locally. One example (from \cite{Ba81a}) is the strongly orientation preserving "angle-doubling" map $y_0:\Omega:=B_1(0)\to B_1(0)$ for $d=2$, defined in complex polar coordinates by $y_0(re^{i\varphi})=re^{2i\varphi}$. 
It satisfies $\deg(y_0;A;z)=2$ for all $z$ close to $0=y_0(0)$ and all open $A\subset \Omega$ with $0\in A$. 
Therefore,  (i) does not hold. 
Moreover, for each $z\in B_1(0)\setminus\{0\}$, $y_0^{-1}(\{z\})=\{x_1,x_2\}\subset \Omega$ contains exactly two (antipodal) points, and thus is not connected as asserted in (ii). In addition, we can always choose $U\in \cO$ with smooth boundary
such that $x_1\in U$ and $x_2\in\partial U$. This violates (iii), as $y(x_1)\in \imloc(y_0;U)$, $y(x_2)\in y_0(\partial U)$
and $y(x_1)=y(x_2)=z$.
\end{rem}
\begin{proof}[Proof of Theorem~\ref{thm:sop-deg1}] 
%Arguments similar to the proof of (ii) are also employed in \cite{Ba81a}. 
%However, we now have to work directly with the image of $y$ (or, better, $\imloc(y,\cdot)$), as %comparison with the image of a given homeomorphism as in \cite{Ba81a} is no longer possible. %Besides, for lack of any given additional (Sobolev) regularity of $y$, we have to handle the %possibility that 
%even smooth boundaries might get mapped to sets with non-empty interior. 

{\bf (i)} Let $z\in y(A)\setminus y(\partial A)$. By Lemma~\ref{lem:sop}, $\Deg(y;\tilde{A};z)\geq 1$ for a suitable $\tilde{A}\in \cO$ with $\tilde{A}\subset A$ and $z\notin y(\partial \tilde{A})$.
Since $\Deg(y;\cdot;\cdot)\geq 0$, we infer that
\begin{align}\label{lemsopdeg1-0a}
	1\leq \Deg(y;\tilde{A};z) \leq \Deg(y;A;z) 
\end{align}
by the additivity of the degree. If, in addition, $z\notin y(\partial U)$, we analogously get that
\begin{align}\label{lemsopdeg1-0b}
	\Deg(y;A;z) \leq \Deg(y;U;z) \leq 1.
\end{align}
Now suppose that $y(\partial U)$ has empty interior, and let $S$ be an arbitrary connected component of $y(A)\setminus y(\partial A)$.
In view of \eqref{lemsopdeg1-0a}, it suffices to show that $\Deg(y;A;\cdot) \leq 1$ on $S$.
As $\Deg(y;A;\cdot)\geq 1$ on $y(A)\setminus y(\partial A)$ by \eqref{lemsopdeg1-0a}
and $\Deg(y;A;\cdot)=0$ on $\RR^d\setminus y(\bar A)$, $y(A)\setminus y(\partial A)$ is open by continuity of the degree. Hence, $S$ is open, too,
and $S\setminus y(\partial U)\neq \emptyset$. Consequently, \eqref{lemsopdeg1-0b} implies that $\Deg(y;A;z_0) \leq 1$ for a $z_0\in S$.
Since $\Deg(y;A;\cdot)$ is constant on $S$, again by continuity of the degree, we conclude that $\Deg(y;A;z) \leq 1$ for all $z\in S$.

{\bf (ii)} %"$\mbf{\cR_y(U) \cap y^{-1}(\{z\})}$ {\bf is connected}": 
The proof is indirect. Let $z\in \imloc(y;U)\setminus y(\partial U)$ and suppose that $\cR_y(U)\cap y^{-1}(\{z\})$ is not connected. We therefore have at least two connected components,
say, $C_z^1$ and $C_z^2$, both of which are compact subsets of $U$.
By Lemma~\ref{lem:sop}, there exist disjoint sets $A^1, A^2\in \cO$ with $A^1\cup A^2\subset U$,
$C_z^j\subset A^j$, $z\notin \partial A^j$ and
$\Deg(y;A^j;z)\geq 1$. Since $z\notin y(\partial U\cup \partial A^1 \cup \partial A^2)$, by additivity of the degree and the fact that $\Deg(y;\cdot;\cdot)\geq 0$ ($y$ is orientation preserving), we obtain that $\Deg(y;\cR_y(U);z)\geq \Deg(y;A^1;z)+\Deg(y;A^2;z)\geq 2$. This contradicts our assumption on the degree.

{\bf (iii)} "$\mbf{\imloc(y;U)\cap  y(\partial U)=\emptyset}$": The proof is indirect. Let $x_0\in \partial U$ and suppose that $z_0:=y(x_0)\in \imloc(y;U)$. By definition of $\imloc(y;U)$, there exists $A\in \cO$, $A\subset U$ with $z_0\notin y(\partial A)$ and $\Deg(y;A;z_0)\neq 0$.
In particular, $y(x)=z_0$ for an $x\in A$ and $x\in \cR_y(U)$. 
Let $C$ denote the connected component of $y^{-1}(z_0)$ containing $x$. 
Since $C$ is compact and $C\subset A$,
$C\subset \cR_y(U)$ by definition of the latter set. In particular,
$C$ has positive distance to $\RR^d\setminus \cR_y(U)\supset \partial U$. By Lemma~\ref{lem:sop}, we can find $A_0\in \cO$ with $A_0\subset A$ such that
\begin{align}\label{lemsopdeg1-1}
	\bar{A}_0\subset  \cR_y(U),~~ z_0\notin y(\partial A_0)~~\text{and}~~\Deg(y;A_0;z_0)\geq 1.
\end{align}
Now choose a sequence $(\tilde{x}_k)\subset U$ with $\tilde{x}_k\to x_0\in \partial U$ as $k\to\infty$. 
Let 
\begin{align*}%\label{lemsopdeg1-2a}
	r(k):=\tfrac{1}{2}\dist{\tilde{x}_k}{\partial U}\leq \tfrac{1}{2}\abs{\tilde{x}_k-x_0} \to 0.
\end{align*}
Since $y$ is strictly orientation preserving, there exist sets $A_k\in \cO$ such that for a $\tilde{z}_k\in \RR^d$,
\begin{align}\label{lemsopdeg1-2}
	A_k\subset B_{r(k)}(\tilde{x}_k), ~~ \tilde{z}_k\notin y(\partial A_k) ~~\text{and}~~\Deg(y;A_k;\tilde{z}_k)\geq 1.
\end{align}
By continuity of the degree and the fact that $y(\partial A_k)$ is compact, \eqref{lemsopdeg1-2} even holds
for all $z_k\in V_k$ in place of $\tilde{z}_k$, in an open neighborhood $V_k$ of $\tilde{z}_k$.
Shrinking $V_k$ if necessary, we can also make sure that $V_k\cap y(\partial A_k)=\emptyset$.
Since, by assumption, $y(\partial U)$ has empty interior while $V_k$ is open, there exists
$z_k\in V_k$ with $z_k\notin y(\partial U)$.
In addition, as $\Deg(y;A_k;z_k)\neq 0$, $A_k\cap y^{-1}(z_k)\neq \emptyset$. Consequently, there exists at least one
connected component $C_k$ of $y^{-1}(z_k)$ in $A_k$, and by definition of $\cR_y(U)$, we also have that $C_k\subset \cR_y(U)$.

On the other hand, $\sup_{x\in A_k} \abs{x-x_0}\leq \tfrac{3}{2}\abs{x_0-\tilde{x}_k}\to 0$ by construction. Since $\bar A_0\subset U$ while $x_0\in \partial U$, 
this implies that $A_k\cap A_0=\emptyset$ for all $k$ large enough, and we also infer that
$z_k=y(x_k)\to y(x_0)=z_0$ by continuity of $y$, for arbitrary $x_k\in C_k\subset A_k\cap y^{-1}(\{z_k\}$.
Arguing as before, we obtain that \eqref{lemsopdeg1-1} also holds for 
$z_k$ instead of $z_0$, for all sufficiently large $k$. As a consequence, besides $C_k\subset A_k\cap \cR_y(U)$, $y^{-1}(\{z_k\})$ has
a second connected component $\hat{C}_k$, now contained in $A_0\cap \cR_y(U)$.
Since we made sure that $z_k\notin y(\partial U)$, this contradicts (ii).

"$\mbf{\cR_y(U)\supset U\setminus y^{-1}(y(\partial U))}$": Let
$z\in y(U)\setminus y(\partial U)$, and
let $C_z$ denote an arbitrary connected component of $y^{-1}(\{z\})$. 
Since $z\notin y(\partial U)$, $C_z$ is compact and has a positive distance to $\partial U$. 
As a consequence of Lemma~\ref{lem:sop}, $C_z\subset \cR(y;U)$.

"$\mbf{\cR(y;U)\subset U\setminus y^{-1}(y(\partial U))}$": By definition, $\cR(y;U)\subset U$, and 
we already know that $\imloc(y;U)\cap y(\partial U)=\emptyset$ and $y(\cR(y;U))=\imloc(y;U)$, the latter by Lemma~\ref{lem:orpres-top} (iii).
Hence, 
\[
	\cR(y;U)\subset y^{-1}(\imloc(y;U))\subset y^{-1}(y(U)\setminus y(\partial U))=U\setminus y^{-1}(y(\partial U)).
\]

"$\mbf{\imloc(y;U)=y(U)\setminus y(\partial U)}$": 
Recall that $\imloc(y;U)=y(\cR_y(U))$ by Lemma~\ref{lem:orpres-top} (iii).
Hence, $\imloc(y;U)\subset y(U)$, and we have just proved that $\imloc(y;U)\cap y(\partial U)=\emptyset$.
We also already know that $U\setminus y^{-1}(y(\partial U))\subset \cR(y;U)$,
which implies that $y(U)\setminus y(\partial U)\subset y(\cR_y(U))=\imloc(y;U)$.

"$\mbf{\imloc(y;U)=\imT(y;U)}$": 
By the definition of $\imT$, $\imT(y;U)\subset y(U)\setminus y(\partial U)$. 
We also have that $\imloc(y;U)\subset \imT(y;U)\cup y(\partial U)$, by the definition of $\imloc(y;U)$ and Lemma~\ref{lem:orpres-top} (ii) with $A_2=U$.
As $y(U)\setminus y(\partial U)=\imloc(y;U)$ due to the previous step, we infer both "$\supset$" and "$\subset$".

"$\mbf{y(\partial U)=\partial (\imloc(y;U))}$": Since $y(U)\setminus y(\partial U)=\imloc(y;U)$ is open and $y(\bar U)$ is closed, it suffices to show 
that $y(\partial U)\subset \overline{\imloc(y;U)}$.
By continuity of $y$, this follows from Lemma~\ref{lem:orpres-top} (iii).

{\bf The case $\mbf{U=\Omega}$ for $\mbf{y\in C(\bar\Omega;\RR^d)}$}: For such a $y$, all notions used above are also defined for $U=\Omega$, and the proof works without changes.
In particular, Lemma~\ref{lem:orpres-top} can be applied with $A_2=U=\Omega$.
\end{proof}
\begin{proof}[Proof of Corollary~\ref{cor:sop-deg1}]
The proof is indirect. Let $K\subset \bar{U}$ be compact with $\interior K=\emptyset$ and suppose there exists a non-empty open set $V\subset y(K)$. In particular, $V\setminus y(\partial U)\neq \emptyset$.

Since $y$ is continuous, there exists $x_0\in K\cap U$ and an open neighborhood $A_0$ of $x_0$
such that $y(A_0)\subset V$. Moreover, $\interior y(A_0)\neq \emptyset$ by the solvability and stability properties of the degree, 
exploiting that $y$ is strictly orientation preserving. 
Choose a non-empty open set
\[
	A\subset \tilde{A}_0:=U\cap y^{-1}(\interior y(A_0))\neq \emptyset\quad\text{such that}\quad A\cap K=\emptyset.
\]
This is possible because $\tilde{A}_0$ is open and $K$ is closed with empty interior.
Now let $x\in A$ and $z:=y(x)$. Since $y(A_0)\subset V\subset y(K)$, $y^{-1}(\{z\})$ also contains a point $x_2\in K$. Let $C_1$ and $C_2$ be the connected components of $y^{-1}(\{z\})\cap U$ with $x\in C_1$ and $x_2\in C_2$, respectively.
By Theorem~\ref{thm:sop-deg1} (ii) and (iii), we have one of the following two possibilities:
\begin{itemize}
\item[(a)] $\bar{C}_1\subset U$ and $C_1=C_2$, or
\item[(b)] $\bar{C}_1 \cap \partial U\neq \emptyset$ and $\bar{C}_2 \cap \partial U\neq \emptyset$.
\end{itemize}
In both cases, $\partial A \cap C_1\neq \emptyset$, because $C_1$ is connected, contains $x\in A$ and has to reach another point outside of $A$. 
With the same argument, we even get that $\partial \tilde{A} \cap C_1\neq \emptyset$ for all open $\tilde{A}\subset A$
with $x\in \tilde{A}$. In other words, $y(x)\in y(\partial \tilde{A})$ for all $x\in \tilde{A}\subset A$.
But for strictly orientation preserving and therefore non-degenerate $y$, this is impossible.
\end{proof}
\begin{proof}[Proof of Corollary~\ref{cor:sop-deg1loc}]
Recall that by Remark~\ref{rem:Ry2},
\begin{align}\label{cord1l-1}
	\cR_y(\Omega)=\textstyle\bigcup_{m} \cR_y(\Omega_m)\quad\text{and}\quad
	\textstyle\imloc(y;\Omega)=\bigcup_{m} \imloc(y;\Omega_m).
\end{align}
The assertions (i)--(iii) now follow from Theorem~\ref{thm:sop-deg1} and 
Corollary~\ref{cor:sop-deg1}.
As to (ii) and (iii) also see Remark~\ref{rem:Ry}.
For (iv), we use \eqref{cord1l-1} and Lemma~\ref{lem:orpres-top} (iii) with $U=\Omega_m$ for each $m$.
\end{proof}
We conclude the section with a result that will be useful to exploit extra regularity 
of deformations with finite distortion while avoiding additional assumptions near the boundary.
\begin{lem}\label{lem:orpres-restrict}
Let $\Omega\subset \RR^d$ be open,
and assume that $y\in C(\Omega;\RR^d)$ is strictly orientation preserving.
If $U\subset \Omega$ is open, the restriction
\[
	\hat{y}:=y|_{\Lambda}:\Lambda\to \RR^d, \quad\text{with}~\Lambda:=\cR_y(U), %\Omega\setminus y^{-1}(y(\partial\Omega))
\]
is strictly orientation preserving on $\Lambda$. 
%with respect to $\hat{\cO}:=\mysetr{A\in\cO}{A\subset\subset \Lambda}$.
If, in addition, $y$ is continuous on $\bar U$, then
\[
	\Deg(y;U;z)=\Deg(\hat{y};\Lambda;z)\quad\text{for all $z\notin y(\partial U)\supset y(\partial \Lambda)$}.
\]
\end{lem}
\begin{proof}
The way we defined strictly orientation preserving maps, 
any restriction like $\hat{y}$ just means fewer sets $A$ to test with and thus trivially inherits this property.
Now suppose that $y$ is continuous on $\bar U$. If $z\notin y(\partial U)$,
\[
	\Deg(y;U;z)=\Deg(y;\emptyset;z)+\Deg(y;\Lambda;z)=\Deg(y;\Lambda;z)
\]
by additivity of the degree. Here, notice that
$y(\partial \Lambda)\subset y(\bar{U}\setminus \Lambda)\subset y(\partial U)$, the latter by Lemma~\ref{lem:orpres-top} (iii).
\end{proof}

\section{Global invertibility in $W^{1,p}$} \label{sec:p>d}

Throughout this section, we will impose the global invertibility constraint $y\in \DEGoneloc$, or $y\in \DEGone$ if $y$ is continuous up to the boundary, for all admissible deformations $y\in W^{1,p}_+$. Recall that by Remark~\ref{rem:invcompare}, this assumption can always be replaced by any of the other invertibility constraints of Section~\ref{sec:AI}, including the Ciarlet-Ne\v{c}as condition $y\in \CNC$ and approximate invertibility on the boundary $y\in\AIB$ or $y\in\AIBloc$. In the latter two cases, we have to assume in addition that $\RR^d\setminus \partial\Omega$ has only two connected components to be able to apply Theorem~\ref{thm:deg-of-def}.

For numerical purposes, $\AIB$ and $\AI(\bar\Omega)$ are more accessible than the other constraints (cf.~\cite{KroeVa19a}).

\subsection{Ball's global invertibility revisited}

The invertibility results obtained in \cite{Ba81a} all rely on the assumption that on $\partial\Omega$, the deformation $y$ (denoted $u$ in \cite{Ba81a}) coincides with 
a continuous $u_0:\bar\Omega\to y(\bar\Omega)$ which is injective on $\Omega$. 
By the results of Section~\ref{sec:top}, it actually suffices to assume that $y\in \DEGone$ instead. This leads to the following generalization of \cite[Theorem 1]{Ba81a}: 
%for strictly orientation-preserving deformations with $\Deg(y;\Omega;\cdot)\leq 1$: 
%Here, recall that for orientation-preserving maps, we always have that $\Deg(y;\cdot;\cdot)\geq 0$.
\begin{thm}\label{cor:Ball}
Let $\Omega\subset \RR^d$ be a bounded Lipschitz domain.
Moreover, let $p>d$, let $y\in C(\bar\Omega;\RR^d)\cap W_+^{1,p}(\Omega;\RR^d)\cap \DEGone$. 
Then we have the following:
\begin{myenum}
\item[(i)] $y(\bar\Omega)=\overline{y(\Omega)\setminus y(\partial\Omega)}$ and 
$y(\partial\Omega)=\partial (y(\Omega)\setminus y(\partial\Omega))$.
\item[(ii)] For every measurable $A\subset \Omega$ and every measurable function $f$, 
\[
	\int_A f(y(x))\det\nabla y(x)\,dx = \int_{y(A)} f(z)\,dz,
\]
as long as at least one of the two integrals exists.
\item[(iii)] For almost every $z\in y(\bar\Omega)$,
$y^{-1}(\{z\})$ consists of a single point.
\item[(iv)] 
If $z\in y(\Omega)\setminus y(\partial\Omega)$, then $y^{-1}(\{z\})\subset \bar\Omega$ 
is a connected set contained in $\Omega$. 
\item[(v)] 
$y^{-1}(y(\partial\Omega))\cap \Omega$ has empty interior, and
if $z\in y(\partial\Omega)$, then each of the connected components of $y^{-1}(\{z\})\cap \Omega$ touches $\partial\Omega$.
%\item If for some $q>d$,
%\begin{align}\label{findistL1}
	%\int_\Omega \abs{(\nabla y(x))^{-1}}^q \det \nabla y(x)\,dx<\infty,
%\end{align}
%then $y$ has a continuous right inverse $y^{-1}\in W^{1,q}(\hat{\cB};\RR^d)$
%on $\hat{\cB}=y(\Omega)\setminus y(\partial\Omega)$,
%with $\nabla(y^{-1})(z)=[\nabla y(y^{-1}(z))]^{-1}$ for a.e.~$z\in \hat{\cB}$.
%If, in addition, $\hat\cB$ satisfies an interior cone condition, then
%$y(\Omega)=\hat{\cB}$ and $y:\Omega\to \hat{\cB}$ is a homeomorphism.
%If $\hat\cB$ is even a Lipschitz domain, then $y:\bar\Omega\to y(\bar\Omega)$ is a homeomorphism.
\end{myenum}
\end{thm}
A proof is given at the end of the subsection.
\begin{rem}
Since $\Omega$ is a Lipschitz domain,
$y\in W^{1,p}(\Omega;\RR^d)$ always has a continuous representative in $C(\bar\Omega;\RR^d)$ by embedding,
and $y$ can always be extended to a function in $W^{1,p}$ on a bigger domain. 
%It is implicitly understood to be used when the choice matters.
\end{rem}
\begin{rem}\label{rem:CN}
With $f\equiv 1$ and $A=\Omega$, Theorem~\ref{cor:Ball} (ii) implies the Ciarlet-Ne\v{c}as condition, cf.~Definition~\ref{def:CNC}.
\end{rem}
\begin{rem}[The case $p=d$]\label{rem:pequald} 
Theorem~\ref{cor:Ball} can be extended to the case $p=d$ with minor modifications, 
see Theorem~\ref{thm:BallAIBloc}.
In this case, continuity of $y$ up to the boundary would be an unnatural extra assumption even for smooth domains. However, 
even if we only have that $y\in W_+^{1,d}(\Omega;\RR^d)$,
we can still follow the proof of Theorem~\ref{cor:Ball} in subdomains compactly contained in $\Omega$, exploiting Remark~\ref{rem:Ry2} and the following facts:
Inside $\Omega$,
deformations $y\in W_+^{1,d}$ automatically have a continuous representative
\cite{VoGo76a} (cf.~\cite[Theorem 5.14]{FoGa95B}, \cite{Sve88a}) and satisfy Lusin's condition (N) 
\cite[Corollary 3.13]{MaZie92a} (cf.~\cite[Theorem 5.32]{FoGa95B}, or \cite[Theorem 4.5]{HeKo14B}).
Even an explicit modulus of continuity can be obtained at any given positive distance from the boundary \cite{Re89B} 
(cf.~\cite[Corollary 5.19]{FoGa95B}). In particular, any such $y$ can still be approximated by smooth functions, simultaneously in $W^{1,d}$ and locally uniformly. (The approximants do not necessarily have positive determinant, though.)
For suitable extensions of the change of variables formulas \eqref{degCOV} and \eqref{areaformula} used below, see \cite[Theorem 5.35 and Theorem 5.34]{FoGa95B} (e.g.). 
\end{rem}

By Theorem~\ref{thm:deg-of-def}, Theorem~\ref{cor:Ball} %(combined with Remark~\ref{rem:pequald} for the case $p=d$)
immediately implies the following variant in the class of approximately invertible maps on the boundary in the sense of Definition~\ref{def:AIB}:
\begin{cor}\label{cor:BallAIB}
Let $\Omega\subset \RR^d$ be a bounded bounded Lipschitz domain 
such that $\RR^d\setminus \partial\Omega$ has exactly two connected components.
If $p>d$ and $y\in W_+^{1,p}(\Omega;\RR^d)\cap \AIB$, 
then the assertions (i)-(v) of Theorem~\ref{cor:Ball} hold.
\end{cor}

In view of Remark~\ref{rem:pequald}, we can also easily obtain an extension of Theorem~\ref{cor:Ball} for the case $p=d$ and 
without requiring $\Omega$ to be Lipschitz, with an analogous proof.
The only difference is a weaker description of $y(\Omega)$, taking into account that we can no longer apply 
Theorem~\ref{thm:sop-deg1} or Lemma~\ref{lem:orpres-top} with $U=\Omega$, only 
Corollary~\ref{cor:sop-deg1loc} which is obtained
by approximating $\Omega$ from inside (see also Remark~\ref{rem:Ry2} and Remark~\ref{rem:interface}).
In (v), the set $y(\Omega)\setminus \imloc(y;\Omega)$ now plays the role of $y(\Omega)\cap y(\partial\Omega)$ which is no longer defined.
\begin{thm}\label{thm:BallAIBloc}
Let $\Omega\subset \RR^d$ be open, let $p\geq d$ and suppose that
$y\in W_+^{1,p}(\Omega;\RR^d)\cap \DEGoneloc$, the latter with respect to a regular inner covering $(\Omega_m)_{m\in\NN}$ of $\Omega$ (see Definition~\ref{def:regcover}).
With $\imT$ and $\imloc$ defined as in Lemma~\ref{lem:orpres-top}, we then have the following:
\begin{myenum}
\item[(i)] $\textstyle\bigcup_{m\in \NN} \imT(y;\Omega_{m}) =\imloc(y;\Omega)\subset y(\Omega)\subset \overline{\imloc(y;\Omega)}$.
\item[(ii)] For every measurable $A\subset \Omega$ and every measurable function $f$, 
\[
	\int_A f(y(x))\det\nabla y(x)\,dx = \int_{y(A)} f(z)\,dz,
\]
as long as at least one of the two integrals exists.
\item[(iii)] For almost every $z\in y(\Omega)$,
$y^{-1}(\{z\})$ consists of a single point.
\item[(iv)] 
If $z\in \imloc(y;\Omega)$, then $y^{-1}(\{z\})$ 
is a compact connected set contained in $\Omega$. 
\item[(v)] $\Omega\cap y^{-1}(y(\Omega)\setminus \imloc(y;\Omega))$ has empty interior, and
if $z\in y(\Omega)\setminus \imloc(y;\Omega)$, then each of the connected components of $y^{-1}(\{z\})\cap \Omega$ touches $\partial\Omega$.
\end{myenum}
\end{thm}

The following lemma links the analytical and topological notions of strictly orientation preserving maps.
\begin{lem}\label{lem:orpres-W1p}
Let $\Omega\subset\RR^d$ be open and bounded, let $p\geq d$ and
let %$y\in C(\Omega;\RR^d)\cap W_{+,\loc}^{1,p}(\Omega;\RR^d)$
$y\in W_{+,\loc}^{1,p}(\Omega;\RR^d)$. Then $y$ is strictly orientation preserving in the sense of Definition~\ref{def:op}.
\end{lem}
\begin{proof}
Let $A\subset\subset \Omega$ be open. As remarked in \cite{Ba81a}, if $z\notin y(\partial A)$
and $V_0$ denotes the connected component of 
$\RR^d\setminus y(\partial A)$ containing $z$,
then the degree can 
be represented as
\begin{align}\label{degcalc-degrep}
	\Deg(y;A;z)=\int_A h(y(x))\,\det\nabla y(x)\,dx,
\end{align}
for any continuous $h:\RR^d\to [0,\infty)$ compactly supported in $V_0$ and with $\int_{V_0} h(z) \,dz=1$. 
Notice that \eqref{degcalc-degrep} is actually a special case of \eqref{degCOV} below which uses that the degree is locally constant. It
thus also extends to the case $p=d$, cf.~Remark~\ref{rem:pequald}.
As an immediate consequence of \eqref{degcalc-degrep},
\begin{align*} 
	\Deg(y;A;z)\geq 0 \text{ for every $z\in \RR^d \setminus y(\partial A)$}.
\end{align*}
Hence, $y$ is orientation preserving in the sense of of Definition~\ref{def:op}.
%(By approximating $\Omega$ from the inside with slightly smaller sets, this also works for $p=d$, cf.~Remark~\ref{rem:pequald}.)
To prove that it is so strictly, first notice that we may assume w.l.o.g.~that $\cL^N(\partial A)=0$, moving to a slightly smaller but more regular set if necessary. For $z\in y(A)\setminus y(\partial A)$, 
we can always choose a function $h$ such that $h(z)>0$. By continuity of $y$, $h\circ y>0$ on a neigborhood of $U\cap y^{-1}(\{z\})$
which automatically has positive measure. Using \eqref{degcalc-degrep} once more, we infer that
\begin{align*} 
	\Deg(y;U;z)>0 \text{ for every $z\in y(A) \setminus y(\partial A)$}.
\end{align*}
In addition, $\cL^d(A)>0$ implies that $\cL^d(y(A))>0$ by the area formula \eqref{areaformula} (see also Remark~\ref{rem:pequald} for $p=d$),
while Lusin's property (N) and $\cL^N(\partial A)=0$ imply that $\cL^d(y(\partial A))=0$.
As long as $A\neq\emptyset$, we know that $\cL^d(A)>0$ since $A$ is open, and consequently, $y(A) \setminus y(\partial A)\neq \emptyset$.
We conclude that $y$ is also strictly orientation preserving in the sense of of Definition~\ref{def:op}.
\end{proof}
\begin{proof}[Proof of Theorem~\ref{cor:Ball}]
%It is actually possible to follow the original proof of \cite{Ba81a}, observing that the homeomorphic extension is essentially only used to show that $\Deg(y;\Omega;\cdot)\leq 1$. 
%Instead, we here give a full 
The topological assertions (i), (iv) and (v) follow from
Corollary~\ref{cor:sop-deg1loc}. Its assumptions hold as described in Remark~\ref{rem:interface}:
$y\in \DEGoneloc$ by Remark~\ref{rem:invcompare} (c), and $y$ is strictly orientation preserving in the sense of Definition~\ref{def:op} by Lemma~\ref{lem:orpres-W1p} below. Due to Lusin's condition (N) \cite{MaMi73a} and the fact that $\cL^{d}(\partial\Omega_m)=0$, we also know that for all $m$, $y(\partial\Omega_m)$ has measure zero and thus empty interior. 

%(and Lemma~\ref{lem:orpres-top}).
%As $\Omega$ is Lipschitz, $\cL^d(\partial\Omega)=0$ and extension of $y$ to a bigger domain is %possible. Hence, Lusins property (N) also applies with $\partial\Omega$, yielding that $\cL^d(y(\partial\Omega))=0$.
In addition, we exploit some change-of-variables formulas 
to show (ii) and (iii) as in \cite{Ba81a}. 
Assertion (ii) is a consequence of a more general change-of-variables formula valid for $y\in W^{1,p}(\Omega;\RR^d)$, $p>d$,
$A\subset\subset \Omega$ open with $\cL^d(\partial A)=0$ and $f\in L^\infty(\RR^d)$:
\begin{align}\label{degCOV}
	\int_A f(y(x))\det\nabla y(x)\,dx = \int_{y(A)} f(z)\Deg(y;A;z)\,dz,
\end{align}
see \cite[Theorem 5.31]{FoGa95a} (e.g.). In our case, for every $z\in y(\Omega)$ and every open $A\subset\subset \Omega$, 
$\Deg(y;A;\cdot)=1$ on $y(A)\setminus y(\partial A)$ by Theorem~\ref{thm:sop-deg1} (i),
and $y(\partial A)$ is always a set of measure zero. 
Hence, \eqref{degCOV} implies assertion (ii) for every open $A\subset\subset \Omega$ with $\cL^d(\partial A)=0$ and $f\in L^\infty(\RR^d)$.
The general case follows with an approximation argument. With the help of the area formula \cite{MaMi73a},
%(see also \cite{MaSwaZie03a}),
\begin{align} \label{areaformula}
	\int_A \abs{\det\nabla y(x)}\,dx = \int_{y(A)} \#y^{-1}(\{z\})\,dz,
\end{align}
where $\#$ denotes the counting measure, (ii) with $f\equiv 1$ implies (iii).
\end{proof}

\subsection{Improved invertibility exploiting finite distortion}

Any map $y$ in 
$W^{1,p}_+(\Omega;\RR^d)$, or, more generally, in $W^{1,p}_{+,\loc}(\Omega;\RR^d)$
is automatically a map of (almost everywhere) \emph{finite distortion}, with \emph{outer distortion} 
\[
	K^O_y(x):=\abs{\nabla y(x)}^d (\det \nabla y(x))^{-1}. 
\]
The \emph{inner distortion} of $y$ is defined as
\[
	K^I_y(x):=\abs{(\nabla y(x))^{-1}}^d \det \nabla y(x)=\abs{\cof \nabla y(x)}^{d} (\det \nabla y(x))^{1-d}. 
\]
Here, $\cof \nabla y$ denotes the cofactor matrix of $\nabla y$, cf.~Remark~\ref{rem:polyconvex}. 
We always have that $(K^O_y)^{d-1}\geq c K^I_y$ with a constant $c=c(d)>0$, because
$\abs{F}^{d-1}\geq c \abs{\cof F}$ for all $F\in\RR^{d\times d}$.

The investigation of maps with finite inner or outer distortion 
was initiated by \cite[Theorem 2]{Ba81a} (for $p>d$) and strongly influenced by \cite{Sve88a} (in particular for $p=d$). 
Their theory is now well developed \cite{HeKo14B}.

If $y\in W^{1,d}_+(\Omega;\RR^d)$, $K^O_y\in L^{q(d-1)}_{\loc}$ with $q>1$ and $y$ is not constant, then it is \emph{open} and \emph{discrete} \cite{ViMa98a}, i.e., $y$ maps open sets to open sets in $\RR^d$ and for any $z\in \RR^d$, $y^{-1}(\{z\})$ does not have accumulation points in $\Omega$. 
A slightly stronger version of the same result was obtained in \cite{Raj11a},
for $K^O_y\in L^{d-1}_{\loc}$ and $K^I_y\in L^{q}_{\loc}$ with $q>1$.

This can be combined with Theorem~\ref{thm:sop-deg1} to 
generalize \cite[Theorem 3.27]{HeKo14B} and the result sketched in \cite[Remark 7.6]{HeKo14B}:
\begin{thm}\label{thm:finoutdistort}
Let $\Omega\subset \RR^d$ be open,  let $p\geq d$ and let $y\in W^{1,p}_{+,\loc}(\Omega;\RR^d)\cap C(\Omega;\RR^d)\cap \DEGoneloc$ such that for a $q>1$,
either $K^O_y\in L^{q(d-1)}_{\loc}(\Omega)$, or
$K^O_y\in L^{d-1}_{\loc}(\Omega)$ and $K^I_y\in L^{q}_{\loc}(\Omega)$. 
Then 
\[
	\text{$y:\Omega\to y(\Omega)$ is a homeomorhpism}
\]
with inverse $y^{-1}\in W^{1,d}_{\loc}(y(\Omega);\RR^d)$. 
If $K^O_y\in L^{d-1}(\Omega)$, we also have that $\nabla y^{-1}\in L^{d}(y(\Omega);\RR^{d\times d})$.

If, in addition, $\Omega$ is bounded, $y\in C(\bar\Omega;\RR^d)$ and $y(\partial\Omega)$ has empty interior, for example if $y\in W^{1,p}_{+}(\Omega;\RR^d)$ with $p>d$ and $\Omega$ is a Lipschitz domain, then $\DEGoneloc$ can be replaced by $\DEGone$ above.
\end{thm}
\begin{rem}\label{rem:GKMS19}
As pointed out in \cite[Section 3]{GraKruMaiSte19a}
for $d=3$ (and easily extended to any $d$), 
openness combined with invertibility almost everywhere 
implies invertibility everywhere. This is essentially equivalent to 
Theorem~\ref{thm:finoutdistort}, although its proof is completely different. For comparison, recall that for $y$ in $W^{1,p}_+$ with $p\geq d$, injectivity a.e.~is equivalent to the Ciarlet-Ne\v{c}as condition $y\in \CNC$. We know that
$\CNC=\DEGoneloc$ by Remark~\ref{rem:invcompare}.
\end{rem}
\begin{proof}[Proof of Theorem~\ref{thm:finoutdistort}]
Let $(\Omega_m)_{m\in \NN}$ be the regular inner covering of $\Omega$ from the definition of $\DEGoneloc$.
Recall that $y$ always has a continuous representative in $\Omega$ and satisfies Lusin's condition (N) (\cite{MaMi73a}; see also Remark~\ref{rem:pequald} if $p=d$).
As a first consequence, $y(\partial \Omega_m)$ has measure zero and thus empty interior for each $m$.
By Lemma~\ref{lem:orpres-W1p}, $y$ is strictly orientation preserving in the sense of Definition~\ref{def:op}. In particular, $y$ is not constant. 
%and
%by Theorem~\ref{thm:sop-deg1} (ii) and (iii) (with $U=\Omega_m$ for arbitrary $m$, cf.~Remark~\ref{rem:interface}), $\imloc(y;\Omega_m)=y(\Omega_m)=y(\cR_y(\Omega_m))$. 
Now let $z\in y(\Omega)$, whence $z\in y(\Omega_m)$ for big enough $m$. Since $y$ is discrete by either \cite[Theorem 1]{ViMa98a} or \cite[Theorem 1]{Raj11a},
all connected components of $y^{-1}(\{z\})$ consist of a single point. By Remark~\ref{rem:Ry}, this entails that $R_y(\Omega)=\Omega$, and by Theorem~\ref{thm:sop-deg1} (ii) (applied with $U=\Omega_m$ for arbitrary $m$; see also Remark~\ref{rem:Ry2}), we infer that $y^{-1}(\{z\})$ is a singleton. 
Hence, $y:\Omega\to y(\Omega)\subset \RR^d$ is bijective, and since it is also continuous, it is a homeomorphism by Theorem~\ref{thm:injcont}.
The Sobolev regularity of $y^{-1}$ 
was shown in \cite[Theorem 5.9]{HeKo14B}.
\end{proof}
If we only control the inner distortion of $y\in W^{1,d}_+(\Omega;\RR^d)$, 
similar results can be obtained given that $y$ is \emph{quasi-light}, i.e., $y^{-1}(\{z\})$ is a compact subset of $\Omega$ for all $z\in y(\Omega)$. Then, $K^I_y\in L^1$ implies that $y$ is either constant or open and discrete \cite{Raj10a}.
Notice that since $y$ is continuous, quasi-light just means that for $z\in y(\Omega)$, there are no connected components of $y^{-1}(\{z\})$ touching $\partial\Omega$. In other words, no continuum connected to the boundary is compressed to a point. This holds by construction if we replace $\Omega$ by the reduced domain $\cR_y(\Omega)$ of Lemma~\ref{lem:orpres-top} (see also Remark~\ref{rem:Ry}). As a result, we can
generalize \cite[Theorem 7.5]{HeKo14B} as follows:
\begin{thm}\label{thm:fininnerdistort}
Let $\Omega\subset \RR^d$ be open, let $p\geq d$, let $y\in W^{1,p}_{+,\loc}(\Omega;\RR^d)\cap C(\Omega;\RR^d)\cap \DEGoneloc$ with $K^I_y\in L^1_\loc(\Omega)$. Then 
\[
	\text{$y:\cR_y(\Omega)\to \imloc(y;\Omega)=y(\cR_y(\Omega))$ is a homeomorhpism}
\]
with $y^{-1}\in W^{1,d}_\loc(\imloc(y;\Omega);\RR^d)$. If $K^I_y\in L^1(\Omega)$, we also get that $\nabla y^{-1}\in L^d(\imloc(y;\Omega);\RR^{d\times d})$.

If, in addition, $\Omega$ is bounded, $y\in C(\bar\Omega;\RR^d)$ and $y(\partial\Omega)$ has empty interior, for example if $y\in W^{1,p}_{+}(\Omega;\RR^d)$ with $p>d$ and $\Omega$ is a Lipschitz domain, then $\DEGoneloc$ can be replaced by $\DEGone$ above and we also 
know that $\imloc(y;\Omega)=y(\Omega)\setminus y(\partial\Omega)=\imT(y;\Omega)$ and $\cR_y(\Omega)=\Omega\setminus y^{-1}(y(\partial\Omega))$.
\end{thm}
\begin{rem}
It is indeed possible that $\Omega\setminus \cR_y(\Omega)\neq \emptyset$ \cite{HeRa13a}.
\end{rem}
\begin{proof}[Proof of Theorem~\ref{thm:fininnerdistort}]
Recall that $\cR_y(\Omega)\subset \Omega$ is open 
and
$y|_{\cR_y(\Omega)}$ is quasi-light by definition of $\cR_y(\Omega)$ (cf.~Remark~\ref{rem:Ry}), 
that $y|_{\cR_y(\Omega)}$ is also strictly orientation preserving in the topological sense like $y$ (Lemma~\ref{lem:orpres-restrict}) and that
$\imloc(y;\Omega)=y(\cR_y(\Omega))$ by Theorem~\ref{thm:sop-deg1} (ii) and Remark~\ref{rem:Ry2}.
In view of these facts, using \cite[Theorem 1]{Raj11a} on $\tilde{\Omega}:=\cR_y(\Omega)$ instead of \cite{ViMa98a} on $\Omega$, 
we argue analogously to the proof of Theorem~\ref{thm:finoutdistort} and infer that $y:\cR_y(\Omega)\to \imloc(y;\Omega)$ is a homeomorphism. 
For the Sobolev regularity of its inverse, first notice that formally, 
$K^I_y(x)=\abs{(\nabla y^{-1})(y(x))}^d \det \nabla y(x)$. 
By \cite[Theorem 5.2]{HeKo14B}, 
%an approximation argument going back to \cite[Theorem 2]{Ba81a} 
we rigorously obtain that $y^{-1}\in W^{1,d}_\loc(\imloc(y;\Omega);\RR^d)$, and $\normn{K^I_y}_{L^1(\Omega)}=\norm{\nabla y^{-1}}_{L^d(\imloc(y;\Omega);\RR^{d\times d})}$ by change of variables.

If $p>d$, $y\in W^{1,p}_{+}(\Omega;\RR^d)$ and $\Omega$ is Lipschitz, then $y$ satisfies Lusin's property (N) and $\cL^d (y(\partial\Omega))=0$.
Lemma~\ref{lem:orpres-top} (iii) and Theorem~\ref{thm:sop-deg1} (iii) with $A=U=\Omega$ provide the asserted properties of $\imloc(y;\Omega)$ and $\cR_y(\Omega)$. 
\end{proof}

\begin{rem}[Connection to light maps] Unlike the argument of \cite{GraKruMaiSte19a}, the proofs of Theorem~\ref{thm:finoutdistort} and Theorem~\ref{thm:fininnerdistort} do not really use that $y$ is open (or discrete), they use that 
\begin{align}\label{totallydc}
	\text{$y^{-1}(\{z\})$ is totally disconnected for each $z\in \RR^d$.}
\end{align}
This is obviously weaker than discreteness. In our setting for strictly orientation preserving maps in $\DEGone$ or $\DEGoneloc$, it also implies openness: \eqref{totallydc} implies that
$y$ is a homeomorphism by the argument in the proof of Theorem~\ref{thm:finoutdistort}.
Maps satisfying \eqref{totallydc} are called \emph{light}. 
%For functions of integrable distortion, this is often shown locally as a step of the proof, in particular in \cite{Raj11a}.
By the Titus-Young theorem \cite[Theorem A]{TiYou62a}, strictly orientation preserving, light maps are always local (but not necessarily global) homeomorphisms on a dense open subset of the domain. 
Connections to our topological results are explained in Remark~\ref{rem:TiYo}.
%A sufficient condition for \eqref{totallydc} is that
%\begin{align}\label{totallydc2}
	%\text{$\cH^1(y^{-1}(\{z\}))=0$ for every $z\in \RR^d$,}
%\end{align}
%where $\cH^1$ denotes the $1$-dimensional Hausdorff measure, or slightly weaker, 
%\begin{align}\label{totallydc3}
	%\text{$\cL^1(e\cdot y^{-1}(\{z\}))=0$ for every $z,e\in \RR^d$ with $\abs{e}=1$,}
%\end{align}
%i.e., the image of $y^{-1}(\{z\})$ in $\RR$ under arbitrary linear projections (realized by $e$) %has measure zero.
%By \cite[Theorem 3.18]{HeKo14B}, \eqref{totallydc2} implies that $y$ is open and discrete for maps %in $W^{1,d}_+$ 
%(even with slightly less integrability for $\nabla y$, namely, $\nabla y\in L^n (\log L)^{-1}$).
\end{rem}
\begin{rem}
For maps that are local homeomorphisms everywhere, also in a suitable neighborhood of each boundary point, invertibility on the boundary is known to imply global invertibility \cite{Wei85a}. This is still true if there are at most finitely many exceptional points inside the domain (and none on the boundary) around which local injectivity does not hold \cite[Theorem 2]{MeiOl63a}. 
The result of \cite{Wei85a} does not require any topologocial restriction on the domain like we do in Theorem~\ref{thm:deg-of-def}. %, and for \cite{MeiOl63a}, it suffices if the boundary is connected.
%For the special case of piecewise affine maps on a simplicial mesh also see \cite{Li14a}.
\end{rem}

\subsection{Existence of homeomorphic minimizers\label{subsec:min}}

In this section, we present an existence result to demonstrate a typical application of our results. 
No attempt is made to achieve maximally general assumptions, and many other variants would be possible, too. In particular,
there are natural generalizations for $d=2$ and $d\geq 4$ instead of $d=3$. 
The special case of controlled outer distortion with deformations
subject to the Ciarlet-Ne\v{c}as condition, (ii) below with $\cG=\CNC$, 
is essentially already covered by the results of \cite[Section 3]{GraKruMaiSte19a}.
Combined with Remark~\ref{rem:invcompare}, their approach can also be used for other 
global invertibility constraints.

The model in the following theorem describes a nonlinearly elastic solid with reference configuration $\Omega$, enclosed in a 
rigid box whose interior is given by $\Lambda$. 
Interpenetration of matter is prevented both locally and globally,
and depending on the shapes of $\Omega$ and $\Lambda$ -- possibly very rough sets -- this can lead to quite interesting, strongly deformed optimal configurations including 
self-contact of the elastic material. All contact is friction-less,
but on a large scale, effectively friction-like forces can still be caused 
if $\Omega$ and $\Lambda$ are rough on comparatively small scales.
\begin{thm}\label{thm:inabox}
Let $p\geq d=3$, let $\Omega,\Lambda\subset \RR^d$ be open and bounded 
and suppose we have a functional (the sum of elastic and potential energy)
\[
	E(y):=\int_\Omega \big(W(\nabla y)+g\cdot y\big) \,dx,\quad E:\cY\to (-\infty,+\infty],
\]
where the class of admissible deformations $y$ is given by
\[
	\cY:=W^{1,p}_+(\Omega;\RR^3)\cap \cG \cap \mysetr{y}{y(\Omega)\subset \bar\Lambda}.
\]
Here, $\cG$ is one of the sets in $\{\DEGoneloc,\CNC,\INV,\AIBloc,\AIloc(\Omega)\}$ (a global constraint preventing self-interpenetration).
If $\cG=\AIBloc$, we also require that $\RR^3\setminus \partial\Omega$ has only two connected components. 
As to the integrand of $E$, we assume that $g\in L^1(\Omega;\RR^3)$ (an external body force, say, gravity)
and the following properties of $W:\RR^{3\times 3}\to [0,+\infty]$ (the stored energy density of the elastic body):
\begin{align}
& \text{$W$ is continuous on }GL_+(3):=\mysetr{F\in \RR^{3\times 3}}{\det F>0}; \label{Wcont}\\
& W(F)=+\infty \quad \text{if and only if}~F\in \RR^{3\times 3}\setminus GL_+(3); \label{Wop}\\
&	W(F)\to+\infty\quad\text{as}~\det F\to 0; \label{Wop2} \\
& W(F) \geq \abs{F}^p  \quad \text{for all}~F\in GL_+(3); \label{Wcoerc} \\
& \text{$W$ is polyconvex (cf.~Remark~\ref{rem:polyconvex})}. \label{Wpc}
\end{align}
Then $E$ attains its minimum in $\cY$, and every minimizer $y^*\in \cY$ is a.e.~injective in $\Omega$.
More can be said if the energy controls the inner or the outer distortion:
\begin{enumerate}
\item[(i)]
If we have in addition that for all $F\in GL_+(3)$,
\begin{align}
& W(F) \geq c \frac{\abs{\cof F}^3}{(\det F)^{2}},
\label{WintegrableKI}
\end{align}
with a constant $c>0$,
then every minimizer $y^*\in \cY$ is a homeomorphism on
$\cR_{y^*}(\Omega)\subset \Omega$, the reduced domain of Lemma~\ref{lem:orpres-top},
and $y^*(\cR_{y^*}(\Omega))$ is open in $\RR^d$ and thus a subset of $\Lambda$.
\item[(ii)] If we even have that for all $F\in GL_+(3)$,
\begin{align}
& W(F) \geq c \left(\frac{\abs{F}^{6}}{(\det F)^{2}}
										+\frac{\abs{\cof F}^{3q}}{(\det F)^{2q}} \right), 
	\label{WintegrableKO}
\end{align}
with constants $q>1$ and $c>0$,
then every minimizer $y^*\in \cY$ is a homeomorphism on $\Omega$,
$y^*(\Omega)$ is open and $y^*(\Omega)\subset \Lambda$.
\end{enumerate}
\end{thm}
\begin{rem}
For example, with $r>0$ and $s\geq 1$, 
\[
	W_1(F):=\abs{F}^p+\frac{1}{(\det F)^{r}},~~~W_2(F):=\abs{F}^p+\frac{1}{(\det F)^{r}}+\abs{\cof F}^s
\]
and 
\[
	W_3(F):=\abs{F}^p+\frac{1}{(\det F)^{r}}+\abs{\cof F}^s+\frac{\abs{F}^6}{(\det F)^{2}}
\]
all satisfy \eqref{Wcont}--\eqref{Wpc}, 
where we set $W_i(F):=+\infty$ if $\det F\leq 0$. 
%The first term of the right hand side of \eqref{WintegrableKO} correspondig to the outer %distortion is polyconvex, and can also be used directly as a summand of $W$, also with higher powers.
They are also frame indifferent in the sense that
$W_i(QF)=W_i(F)$ for all $F\in GL_+(3)$ and all rotations $Q\in SO(3)$.
Moreover, \eqref{WintegrableKI} holds
\begin{itemize}
\item in case $W=W_1$, if $p>6$ and $r\geq \frac{2p}{p-6}$, and
\item in case $W=W_2$ or $W=W_3$ for any $p\geq 3$, if $r>2$ and $s\geq \frac{3r}{r-2}$,
\end{itemize}
by Young's inequality. For $W=W_1$ and $W=W_3$,
strict inequalities yield \eqref{WintegrableKO}.
%For the first case, we also used the trivial estimate $\abs{\cof F}^3\leq C \abs{F}^6$ with a constant $C>0$.
\end{rem}
\begin{rem}
If $\Omega\setminus \cR_{y^*}(\Omega)\neq \emptyset$, then 
by the results of Section~\ref{sec:top},
it consists of connected sets $C\subset \Omega$
in such a way that for each $C$, $C$ touches $\partial \Omega$ and $y^*$ compresses $C$ to a point 
in $y^*(\Omega)\setminus y^*(\cR_{y^*}(\Omega))$ ($\subset y^*(\partial\Omega)$, if $y^*$ happens to be continuous up to the boundary). 
We also know that $\Omega\setminus \cR_{y^*}(\Omega)$ is small in the sense that it has empty interior; 
it even has measure zero, because $y^*$ is a.e.~invertible (because any $y\in \cY$ automatically satisfies the Ciarlet-Ne\v{c}as condition $y\in \CNC$, see Remark~\ref{rem:invcompare}). For more information on the set $\cR_{y}(\Omega)$ see Lemma~\ref{lem:orpres-top},
Theorem~\ref{thm:sop-deg1} and the Remarks~\ref{rem:Ry} and~\ref{rem:Ry2}.
\end{rem}
\begin{rem}
The theorem is written mostly for the case $p=d$; if we assumed and exploited $p>d$ and a Lipschitz domain $\Omega$, the results would be slightly stronger and simpler to state, and we could also admit the constraints $\cG=\AIB$ (if $\Omega$ has connected boundary), $\cG=\AI(\bar\Omega)$ or $\cG=\DEGone$.
\end{rem}
\begin{rem}\label{rem:polyconvex}
Polyconvexity as required in \eqref{Wpc} means that there exists a function
\[
\bald
	&h:\RR^{3\times 3}\times \RR^{3\times 3}\times (0,\infty)\to \RR~~\text{convex, such that} \\
	&W(F)=h(F,\cof F,\det F)~~\text{for all}~F\in GL_+(3),
\eald
\]
cf.~\cite{Ba77a}. Here, $\cof F$ denotes the cofactor matrix, i.e., a matrix formed of the determinants of all $2\times 2$-submatrices of $F$. 
Usually, they are ordered and given suitable signs so that $(\cof F)^T F=\det F$, but this is irrelevant for our purposes.
\end{rem}
\begin{proof}[Proof of Theorem~\ref{thm:inabox}]
The existence of minimizers is a standard application of the direct method. 
First observe that there always is $y\in\cY$ with $E(y)<\infty$, for instance
$\hat{y}(x):=z_0+\lambda x$, where $z_0\in \Lambda$ is chosen arbitrarily but fixed 
and $\lambda>0$ is small enough so that $\hat{y}(\Omega)\subset \Lambda$, exploiting that 
$\dist{z_0}{\partial\Lambda}>0$. 

By the constraint $y(\Omega)\subset \bar\Lambda$,
$\cY$ is bounded in $L^\infty$. Hence, the linear force term 
$y\mapsto \int_\Omega g\cdot y\,dx$ is well defined, and by dominated convergence,
it is also sequentially continuous in $\cY$, first with respect to pointwise convergence almost  everywhere and then also with respect to weak convergence in $W^{1,p}$.
 
Using \eqref{Wcoerc} and the a-priori bound in $L^\infty$, 
it is not difficult to show the coercivity estimate 
\begin{align}\label{thmbox-1}
	E(y)\geq c_1 \norm{y}_{W^{1,p}}^p-c_2 \quad\text{for $y\in W^{1,p}(\Omega;\RR^d)$ with $y(\Omega)\subset \bar\Lambda$,}
\end{align}
where $c_1>0$ and $c_2$ are real constants. By arguments of \cite{Mue90a} (or \cite{Ba77a} if $p>d=3$), 
\begin{align} \label{thmbox-2}
\bald
	\text{$y\mapsto \textstyle\int_\Omega W(\nabla y)\,dx$ is weakly sequ.~lower semicontinuous}&\\
	\text{on $W^{1,p}(\Omega;\RR^d)\cap \{E<\infty\}$}&.
\eald
\end{align}
Here, the essential ingredients for the proof of \eqref{thmbox-2} 
are the weak continuity of $y\mapsto \det\nabla y$, $W^{1,p}_+\to L^1$ on compact subsets of $\Omega$,
the weak continuity of $y\mapsto \cof \nabla y$, $W^{1,p}_+\to L^{p/2}$, and the 
convexity of $h$ (the polyconvexity of $W$, cf.~Remark~\ref{rem:polyconvex}).

As a consequence of \eqref{thmbox-1}, any sequence $(y_k)\subset \cY$ with $E(y_k)\to I:=\inf_{y\in\cY}E(y)<\infty$
is bounded and has a weakly convergent subsequence in $W^{1,p}$, say, $y_k\rightharpoonup y^*$.
Due to \eqref{thmbox-2}, $E(y^*)\leq \lim_k E(y_k)=I$. 
As we also have that $y_k\to y^*$ locally uniformly (by embedding if $p>d=3$, or by Remark~\ref{rem:pequald} if $p=3$), 
we obtain that $y^*(\Omega)\subset \bar\Lambda$. In addition, $\det\nabla y^*>0$ a.e.~by \eqref{Wop}, because $E(y^*)<\infty$.
Hence, $y^*\in W^{1,p}_+$. Finally, 
$y^*\in \cG$, because $\cG\cap W^{1,p}_+$ is weakly sequentially closed in all cases (see Section~\ref{sec:AI}).
Altogether, $y^*\in\cY$ is a minimizer.

The remaining assertions (and some additional properties) follow from Theorem~\ref{thm:fininnerdistort} and Theorem~\ref{thm:finoutdistort}, repectively.
The assumption on $K^I_{y^*}$ or $K^O_{y^*}$ are obtained from 
\eqref{WintegrableKI} or \eqref{WintegrableKO}, and by Remark~\ref{rem:invcompare}, the theorems can be applied for all possible choices of $\cG$.
\end{proof}

\appendix

\section{The problem of homeomorphic extension} \label{sec:A:homext}

When working with injective continuous maps, it is good to keep in mind the following two well-known facts.
\begin{lem}\label{lem:injcont} 
Let $X,Z$ be topoogical spaces and $y:X\to Z$ continuous and injective, and suppose that $X$ is compact.
Then $y:X\to y(X)$ is a homeomorphism, where $y(X)$ is endowed with the trace topology of $Z$.
\end{lem}
\begin{proof}
By continuity of $y$, $y(X)$ is also compact. 
Open sets in $X$ and $y(X)$, respectively, are thus exactly the complements of compact sets.
Since $y:X\to y(X)$ is bijective and maps compact sets to compact sets, it therefore
also maps open sets in $X$ to open sets in $y(X)$.
\end{proof}
The statement above does \emph{not} mean that $y$ maps opens sets in $X$ to open sets in $Z$,
because $y(X)$ is usually not open in $Z$. In $\RR^d$, more can be said:
\begin{thm}\label{thm:injcont} 
Let $\Omega\subset\RR^d$ be open, $\Omega\neq \emptyset$, and let $y:\Omega\to \RR^n$ be injective and continuous. Then $n\geq d$. Moreover, $y(\Omega)$ is open in $\RR^n$ if and only if $n=d$. For $n=d$, $y:\Omega\to y(\Omega)$ is a homeomorphism.
\end{thm}
\begin{proof}
This is a combination of the openness (invariance of domain) and invariance of dimension theorems based on the topological degree, see \cite[Thm.~3.30, Cor.~3.31 and Cor.~3.32]{FoGa95B} (e.g.).
The last assertion is a consequence of the others which also hold for arbitrary open subsets of $\Omega$, thereby proving that $y$ maps open subsets of $\Omega$ to open sets in $\RR^n=\RR^d$. 
\end{proof}

\subsection{Homeomorphic extension versus Schoenflies extension}

We are here mainly interested in homeomorphic extension
for functions given on the boundary of some domain in $\RR^d$. 
In the most general form, this problem reads as follows:
\begin{problem}[Homeomorphic extension problem in $\RR^d$]\label{prob:ext0}
Let $\Omega\subset \RR^d$ a bounded domain, and suppose that
$y: \partial\Omega\to \RR^d$ is continuous and injective.
Is there a homeomorphism $h:\bar\Omega \to h(\bar\Omega)\subset \RR^d$ such that $h=y$ on $\partial\Omega$?
\end{problem}
Additional assumptions on topological nature of $\Omega$ and $\partial\Omega$ are typically added, as
it is well-known that simple counterexamples exist when $\Omega$ is topologically complicated. For instance,
on an annulus, homeomorphic extension is impossible if the winding numbers of $y$ 
on the two boundary pieces are not the same (e.g., one clockwise and the other counterclockwise). 

A close relative is the following question which is much more widely studied in the literature:
\begin{problem}[Schoenflies extension problem]\label{prob:ext}
Let $y: S^{d-1}\to S^d$ continuous and injective.
Is there a homeomorphism $h:S^d\to S^d$ such that $h(S^{d-1})=\Sigma^{d-1}:=y(S^{d-1})$, where $S^{d-1}\subset S^d$ is interpreted as the equator 
of the $d$-dimensional sphere $S^d\subset \RR^{d+1}$?
\end{problem}
The Schoenflies extension problem imposes a restriction on the topological type of admissible domains (one of the half-spheres separated by the equator) -- it must be a topological ball -- which is also commonly used for the homeomorphic extension problem. Apart from that, the two problems are essentially equivalent:
\begin{prop}[Schoenflies versus homeomorphic extension]\label{prop:SchoenHom}~\\
Suppose that $\Omega\subset \RR^d$ is a bounded domain such that there exists 
a homeomorphism $\gamma:\bar\Omega\to \bar{H}$, where $H\subset S^{d}\subset \RR^{d+1}$ is one of the two hemispheres of $S^d$ separated by the "equator" $S^{d-1}$.
Moreover, let $y: S^{d-1}\to S^d$ continuous and injective and let $\delta:S^d\to \RR^d\cup\{\infty\}$ be a homeomorphism with $\delta(y(S^{d-1}))\subset \RR^d$, where $\RR^d\cup\{\infty\}$ denotes the one-point compactification of $\RR^d$. 
Then we have the following for $\tilde{y}:=\delta\circ y \circ \gamma$:
\begin{enumerate}
\item[(i)]  
If there exists Schoenflies extension $h:S^d\to S^d$ of $y$ as in Problem~\ref{prob:ext}, 
then a homeomorphic extension $\tilde{h}$ of $\tilde{y}$ as in Problem~\ref{prob:ext0} exists, too.
\item[(ii)] Conversely, if no Schoenflies extension $h:S^d\to S^d$ of $y$ exists, 
a homeomorphic extension also fails to exist
for one of the following two maps:
\begin{enumerate}
\item[(a)] $\tilde{y}:\partial\Omega\to \RR^d$, or
\item[(b)] $\hat{y}:=\tilde{y}\circ \iota:\partial \hat{\Omega}\to\RR^d$, with
$\hat{\Omega}$ denoting the bounded connected component of $\RR^d\setminus \iota(\partial\Omega)$. 
\end{enumerate}
\end{enumerate}
Here, $\iota: \RR^d\cup\{\infty\}\to \RR^d\cup\{\infty\}$, $\iota(x):=\abs{x-x_0}^{-2}(x-x_0)$, 
is the inversion map with respect to a point $x_0\in \RR^d$; for (ii), we choose an arbitrary but fixed $x_0\in \Omega$.
\end{prop}
\begin{proof}
{\bf (i)} $\tilde\varrho:=h^{-1}\circ y:S^{d-1}\to S^{d-1}$ is a homeomorphism of the equator $S^{d-1}$ onto itself.
It has an explicit "radial" homeomorphic extension $\varrho:S^d\to S^d$. 
Using cylindrical coordinates $(x,t)\in S^{d-1}\times [-1,1]$, 
it is given by 
\[
	\varrho\big((1-t^2)^\frac{1}{2}x,t\big):=\big((1-t^2)^{\frac{1}{2}} \tilde\varrho(x),t\big)\in S^d\subset \RR^d\times \RR.
\]
Using the inversion map $\iota$ with respect to a point $x_0$ in the bounded connected component of $\tilde{y}(\partial\Omega)\subset \RR^d$,
we now define $\tilde{h}:=\delta\circ h\circ \varrho \circ \gamma$ or $\tilde{h}:=\iota\circ \delta \circ h \circ \varrho \circ \gamma$.
One of the two options satisfies $\infty\notin \tilde{h}(\bar\Omega)$, and for this choice, $\tilde{h}|_{\bar\Omega}$ is a homeomorphic extension of 
$\tilde{y}$ in the sense of Problem~\ref{prob:ext0}.

{\bf (ii)} If a Schoenflies extension $h:S^d\to S^d$ of $y$ does not exists,
then extension already fails in one of the two hemispheres of $S^d$ separated by $S^{d-1}$, 
either $\bar{H}=\gamma(\bar\Omega)$ or $S^d\setminus H$. (Otherwise, the extensions to the hemispheres can be glued to a Schoenflies extension, 
after first matching their parametrization of $y(S^{d-1})$ using the radial homeomorphic extension of the proof of (i).)
Accordingly, for either $\tilde{y}$ or $\hat{y}=\tilde{y}\circ \iota$, where $\iota$ is taken with respect to an $x_0\in \Omega$, 
there exists no homeomorphism defined on the closure of $\Omega$ or $\hat\Omega$, respectively, who maps the boundary of its domain to 
$\tilde{y}(\partial\Omega)=\hat{y}(\partial \hat\Omega)$. In particular, either $\tilde{y}$ or $\hat{y}$ has no
homeomorphic extension to its domain in the sense of Problem~\ref{prob:ext0}.
\end{proof}

%Compared to Problem~\ref{prob:ext0}, this corresponds to a setting where $\Omega$ is the unit ball, whose boundary is given by $S^{d-1}$, and $y$ maps into $S^d$ instead of $\RR^d$.
%Replacing the $d$-dimensional sphere $S^d$ as the ambient space in the image with $\RR^d$ is not difficult, since $S^d$
%can be identified with the compactification of $\RR^d$ by adding a point at infinity. Moreover, while a Schoenflies extension 
%$h$ is only required to share the image of $y$ on $S^{d-1}$ as a set, but not its parametrization by $y$,
%the latter can always be achieved in a second step: Given a Schoenflies extension, we can reduce Problem~\ref{prob:ext0} with for the unit ball $\Omega=B$ to the much simpler case $y(\partial B)=\partial B$, where the straightforward radial extension to $B$ is available. The proof of Theorem~\ref{thm:homext} below explains this connection more precisely in a particular case.

\subsection{Known results and counterexamples}

For $d=2$, the answer to Problem~\ref{prob:ext} is affirmative, given by the classical Schoenflies Theorem. Extension theorems for more regular classes of invertible functions are also known in this case, for instance bi-Lipschitz extensions \cite{Tu80a,DaPra15a}. Recently, an extension result (also) valid in the class of Sobolev homeomorphisms has been established in \cite[Theorem 4 and Corollary 5]{IwaOnn19a}. This is based on $p$-harmonic extension and even smooth in $\Omega$.

For $d\geq 3$, the situation is significantly more complicated. In general, a Schoenflies extension can fail to exist,
for instance in case of Alexander's horned sphere \cite{Al24a}. 
However, the result can be recovered for $d\geq 3$ if the embedded sphere $\Sigma^{d-1}:=y(S^{d-1})$ is
\emph{locally flat}:
\begin{thm}[Generalized Schoenflies Theorem {\cite[Theorem 4]{Bro62a}}] \label{thm:GenSchoenflies}~\\
Let $\Sigma^{d-1}\subset S^d$ be a homeomorphic embedding of $S^{d-1}$ which is \emph{locally flat}
in the following sense: 
\[
\bald
	&\text{For each $x_0\in \Sigma^{d-1}$, there exists a neighborhood $V$ of $x_0$ in $S^d$ and}\\
	&\text{a homeomorphism $\zeta:V\to \zeta(V)\subset S^d$ s.t. }\text{$\zeta(V\cap \Sigma^{d-1})\subset S^{d-1}$,}
\eald
\]
where $S^{d-1}$ is interpreted as the equator of $S^{d}$.
Then $\Sigma^{d-1}$ is \emph{flat}, i.e., there is a homeomorphism $h:S^d\to S^d$ such that
$h(\Sigma^{d-1})=S^{d-1}$. 
\end{thm}
There are also variants of the Generalized Schoenflies Theorem that require higher regularity of $\Sigma^{d-1}$ instead of assuming a locally flat embedding.
In particular, this is possible for the piecewise affine (polyhedral) \cite{Al24b} or diffeomorphic \cite{Ma61a} case. 
%As pointed out in \cite{Bro62a}, these assumptions can also be used to prove local flatness.  
As pointed out in \cite[Example 3.10 (5)]{LuuVai77a} for $d=3$, bi-Lipschitz regularity is \emph{not} enough.
\begin{rem}[Homeomorphic extension may fail for $d\geq 3$]\label{rem:Fox-Artin}
In view of Proposition~\ref{prop:SchoenHom}, \cite[Example 3.10 (5)]{LuuVai77a} also entails that for $d=3$, homeomorphic extension is in general impossible even if the given boundary homeomorphism $y:\partial\Omega\to \Sigma^{d-1}:=y(\partial\Omega)$ is bi-Lipschitz.
%In particular, the bi-Sobolev regularity $y\in W^{1,\infty}$ and $y^{-1}\in W^{1,\infty}$ on the boundary is not enough to guarantee the existence 
%of homeomorphic extensions.
\end{rem}
\begin{rem}
The example of \cite{LuuVai77a} is based on the Fox-Artin arc \cite[Example 1.1]{FoArt48a}, a bi-Lipschitz embedding of a compact interval into $\RR^3$ whose image has a complement which is not simply connected. 
By thickening it, surrounding the original interval by a domain consisting of two thin cones back-to-back with tips at the two end points of the interval, the self-similar construction yields a bi-Lipschitz mapping of the domain boundary onto a surface in $\RR^3$. This surface is a topological $2$-sphere, and from the Fox-Artin arc, it inherits that the unbounded component of its complement is not simply connected. In particular, a Schoenflies extension (after identifying $\RR^3\cup \{\infty\}$ with $S^3$) is impossible because its existence would imply that both halves of $S^3$ separated by the surface are topological $3$-balls which are simply connected. 
\end{rem}

On the other hand, if we look for a solution of Problem~\ref{prob:ext} when the embedding of the sphere is known to be locally flat with higher regularity given in the whole neighborhood of its image $\Sigma^{d-1}$, then
this regularity sometimes can be carried over to a suitable extension. 
In particular, this is possible in the bi-Lipschitz case \cite[Theorem 7.7]{LuuVai77a}, or
for the second order bi-Sobolev homeomorphisms
where both the function and its inverse are in $W^{2,p}$ with $1\leq p<d$ \cite{Go11a}. 
As far as I know, there is no comparable result for bi-Sobolev homeomorphisms in $W^{1,p}$ (yet?),
only the theory for maps with finite distortion \cite{HeKo14B} which is conceptually 
closer to regularity theory than to extension results.
%, unless an orientation preserving extension is already given as in \cite[Theorem 7.5]{HeKO14B} (e.g.) for $p=d$,
%However, in this case there is no need to assume any kind of local flatness, as shown in Theorem
In the diffeomorphic category, extensions starting from locally flat embeddings face another obstacle in higher dimensions, the possible existence of exotic spheres, for example for $d=8$ (7-dimensional spheres) \cite[Theorem 3.4]{Mi56a}.

\subsection{Homeomorphic extension for $C^1$ functions on Lipschitz domains}

As shown in Proposition~\ref{prop:SchoenHom}, a solution to Problem~\ref{prob:ext} can be used to build homeomorphic extensions of maps $y|_{\partial\Omega}$, at least if $\Omega$ is homeomorphic to the closed unit ball. A more practical application in the same spirit is given below, using Theorem~\ref{thm:GenSchoenflies} to obtain a homeomorphic extension of a $C^1$-deformation on a Lipschitz domain which is invertible on the boundary. Despite the similarity, it does not directly follow from the Schoenflies extension for a $C^1$ map outlined in \cite[p.11]{Ma61a}, because we would first have to transform the given Lipschitz domain to the unit ball. This is possible, but the transformation is only bi-Lipschitz and we would lose the crucial $C^1$ regularity (cf.~Remark~\ref{rem:Fox-Artin}). 
\begin{thm}\label{thm:homext}
Suppose that $\Omega\subset \RR^d$ is a Lipschitz domain such that $\bar\Omega$ is homeomorphic to the closed unit ball.
Moreover, let $y\in C^1(\bar\Omega;\RR^d)$ such that $y|_{\partial\Omega}$ is injective and
$\det \nabla y\neq 0$ on $\partial\Omega$.
Then $y|_{\partial\Omega}$ admits a homeomorphic extension to $\bar\Omega$.
\end{thm}
\begin{proof}
The proof is based on Theorem~\ref{thm:GenSchoenflies} and Proposition~\ref{prop:SchoenHom} (i). To apply the theorem, we identify $S^d$ with the one-point compactification $\RR^d\cup \{\infty\}$ of $\RR^d$. In this sense, $\RR^d\subset S^d$ (homeomorphically embedded), and $y$ maps $\bar\Omega$ to $\RR^d\subset S^d$,
and $\Sigma^{d-1}:=y(\partial\Omega)$ is a homeomorphic embedding of a topological $(d-1)$-dimensional sphere into $\RR^d\subset S^d$.
To see that this embedding is also locally flat in the sense of Theorem~\ref{thm:GenSchoenflies}, it suffices to
define local bi-Lipschitz extensions of $y\in C^1(\bar\Omega;\RR^d)$ in a neighborhood of each boundary point $x_0\in \partial\Omega$.
Here, notice that $\partial\Omega$ is locally the graph of a Lipschitz function. This implies that $\partial\Omega$ is
locally flat, and we may therefore assume that the homeomorphism mapping $\Omega$ to the unit ball is defined on a whole neighborhood of $\bar\Omega$.

Since $\partial\Omega$ is Lipschitz, we can choose a cylidrical neigborhood of the form 
$C_\eps(x_0):=D_\eps(x_0)+(-\eps,\eps)\nu\subset \RR^d$ with a unit vector  $\nu=\nu(x_0)\in \RR^d$ and a 
$(d-1)$-dimensional disc $D_\eps(x_0)$ of radius $\eps$, centered at $x_0$ and perpendicular to $\nu$. For $\eps>0$ small enough and an appropriate choice of $\nu$,
$\partial\Omega\cap C_\eps(x_0)$ can be represented as the graph of a Lipschitz function $g:D_\eps(x_0)\to (-\eps,\eps)$, 
such that $\Omega\cap C_\eps(x_0)=\mysetr{x'+t\nu}{t<g(x')}$.
We can now extend $y|_{\partial\Omega}$ to a function $\hat{y}:C_\eps(x_0)\to \RR^d$ by setting
\[
	\hat{y}(x'+t\nu):=y(x'+g(x')\nu)+(t-g(x'))Dy(x_0)\nu.
\]
for $x'\in D_\eps(x_0)$ and $t\in (-\eps,\eps)$. Close to $x_0$, this extension 
divides $C_\eps(x_0)$ into surfaces of the form $\partial\Omega+s\nu$ 
(parametrized by $x'+g(x')\in \partial\Omega$ and $s=t-g(x')$) and
maps each such surface onto $y(\partial\Omega)+s Dy(x_0)\nu$, a
shifted copy of $y(\partial\Omega)$. 

To see that $\hat{y}$ is bi-Lipschitz in a neighborhood of $x_0$, the key observation is the following:
Just like $\nu$ and $\partial\Omega$, $Dy(x_0)\nu$ and the surface $y(\partial\Omega)$ always form an angle bounded away from zero as long as we remain close enough to $x_0$, because $Dy$ is continuous and $Dy(x_0)$ is invertible. As an immediate consequence, 
$\partial\Omega\to\RR^d$, $\sigma \mapsto y(\sigma)+s Dy(x_0)$ is injective near $x_0$ for each $s$, and
in a small enough neighborhood $V$ of $y(x_0)$ in $\RR^d$, we also have that
\[
	V\cap [y(\partial\Omega)+s Dy(x_0)\nu]\cap [y(\partial\Omega)+s_2 Dy(x_0)\nu]=\emptyset\quad\text{for}~s_1\neq s_2.
\]
Hence, $\hat{y}$ is injective near $x_0$. Further details are omitted.
%To see that this extension is indeed bi-Lipschitz, suppose for simplicity that $x_0=0$ and $\nu=e_d=(0,\ldots,0,1)$.
%As $g$ is Lipschitz, $Dg(x')[z']\leq L\abs{z'}$ for all $x'\in D_\eps(x_0)$ and $z'\perp \nu$.
%with a constant $\alpha<1$.
%
%since $Dy$ is continuous and $Dy(x_0)$ is invertible, the angle between the surface $x'\mapsto x'+\nu g(x')$ and the direction $D_\eps(x_0)[\nu]$ 
%is bounded away from zero, at least if $\abs{x'-x_0}<\eps$ is small enough:
%\[
%\bald
	%&Dy(x_0)[(x'-x_0+\nu Dg(x_0)[x'-x_0])]\cdot Dy(x_0)[\nu]
	%
	%\\
	%&\leq \abs{Dy(x_0)^{-1}} \abs{(x'-x_0+\nu Dg(x_0)[x'-x_0])\cdot \nu}+o(\abs{x'-x_0})\\
	%&\geq \abs{Dy(x_0)^{-1}} \abs{Dg(x_0)[x'-x_0])\cdot \nu}+o(\abs{x'-x_0})\\
%\eald
%\]
%Here, we used that $y$ is $C^1$ and that $g$ is Lipschitz and thus also in $W^{1,\infty}$.

Theorem~\ref{thm:GenSchoenflies} now gives us a Schoenfliess extension of $\Sigma^{d-1}:=y(\partial\Omega)\subset \RR^d\cup\{\infty\}\cong S^d$,
and by Proposition~\ref{prop:SchoenHom} (i), this implies the existence of a homeomorphic extension of $y|_{\partial\Omega}$ to $\bar\Omega$.
\end{proof}

\section{Counterexamples for domains in with holes}\label{sec:A:topological}
The following examples illustrate that the assumption that $\RR^d\setminus \partial\Omega$ has only two connected components cannot be dropped
in Theorem~\ref{thm:deg-of-def}. 
For simplicity, they are all constructed for $d=2$, but they have straightforward equivalents in higher dimensions, still using 
domains with holes. In particular, it does not really matter whether $\Omega$ is simply connected or not.
In both examples, the explicit values asserted for the degree
are always taken at a suitable regular value of $y$ with just one pre-image $x_0\in \Omega$, and are therefore given
as the sign of $\det\nabla y(x_0)$. Geometric intuition provides a good heuristic, observing whether or not the local deformation is orientation-preserving. 
If yes, the sign is positive, otherwise negative.
\begin{ex}\label{ex:top1}
Take the annulus $\Omega:=B_2(0)\setminus \bar B_1(0)\subset \RR^2$ and for 
$x\in\Omega$ consider $y\in W^{1,\infty}(\Omega;\RR^2)$, 
\[
	y(x):=2\frac{\abs{x}-1}{\abs{x}}x+\frac{2-\abs{x}}{\abs{x}}\left(x+(3,0)\right). 
\]
As defined, $y$ keeps the outer boundary $\partial B_2(0)$ fixed while 
translating the inner boundary onto $y_1(\partial B_1(0))=(3,0)+\partial B_1(0)$.
In particular, $y|_{\partial\Omega}$ is invertible, but it maps $\partial\Omega$ to two circles that lie outside of each other.
Now, $\cB=B_2(0)\cup [(3,0)+B_1(0)]$ and $\sigma$ changes sign; more precisely, $\Deg(y;\Omega;(0,0))=1$ while $\Deg(y;\Omega;(3,0))=-1$.
\end{ex}
It is also not enough to have that $\partial\Omega$ is connected:
\begin{ex}\label{ex:top2}
Take a fixed unit vector $e\subset \RR^2$ and truncated open cones of the form
\[
	\hat{V}(\alpha,r):=O(r,\alpha)\cup\mysetr{x\in \RR^2}{x\cdot e> (1-\alpha)\abs{x},~\abs{x}<r},~~ \alpha,r>0,
\]
where $O(r,\alpha)\subset \RR^s$ denotes the unique open ball 
which touches the surface of the unbounded cone tangentially at $\abs{x}=r$.
Consequently, $V_{\alpha,r}$ has a boundary of class $C^1$ everywhere except at its tip in the origin.
We create a domain by removing a smaller cone from a bigger one sharing the same tip:
\[
	\Omega:=V_2\setminus \bar V_1,\quad \text{where }V_s:=\hat{V}\big(\tfrac{s}{3},s\big).
\]
As a first step, we now consider a map $y\in W^{1,\infty}(\Omega;\RR^2)$ which keeps the outer part of the boundary fixed while flipping the inner part outside,
with affine interpolation on suitable rays in between. For its explicit definition, 
the flip is realized by the reflection $R$ across $\{x\cdot e=0\}$, $R x:=x-2(x\cdot e)e$, and we use the
(nonlinear) projections $Q(x)$ and $P(x)$ onto the inner and the outer boundary, respectively, 
along lines perpendicular to the inner boundary $\partial V_{1}$: For all $x\in\Omega\setminus \{0\}$,
\[
	Q(x)\in \partial V_{1},~~P(x)\in \partial V_{2}~~\text{and}~~
	P(x)-Q(x)\perp \partial V_{1}~~\text{at}~Q(x).
\]
Notice that $Q,P:\Omega \to \RR^2$ are Lipschitz and thus in $W^{1,\infty}$ (even $C^1$ away from $x=0$),
and both converge to the origin as $\abs{x}\to 0$, $x\in\Omega$.
We now can define
\[
\bald
	y(x):=  \frac{\abs{x-Q(x)}}{\abs{P(x)-Q(x)}} P(x) + \frac{\abs{P(x)-x}}{\abs{P(x)-Q(x)}} R[Q(x)] &\\
	= \left(1+\frac{\abs{P(x)-x}}{\abs{P(x)-Q(x)}}\right)\big(R[Q(x)]-P(x)\big)&.
\eald
\]
The latter representation shows that $y$ is Lipschitz also at the origin.

This construction does not yet contradict Theorem~\ref{thm:deg-of-def}, because as a matter of fact, 
$\Deg(y;\Omega;\cdot)=1$ on both components of $\cB(\RR^d\setminus y(\partial\Omega))
=V_2\cup (-V_1)$. 
In any case, with a second deformation that squeezes the line orthogonal to $e$ in the image to the origin and simultaneously reflects the 
half-space $\{x\cdot e<0\}$ containing $-V_1$ across the line in direction $e$, 
we can make the degree change sign at the value $-e\in -V_1$ while keeping it fixed at $e\in V_2$. More precisely,
for 
\[
	\hat y:=F\circ y\quad\text{with}~F(z):=(z\cdot e) e+(z\cdot e)(z\cdot e^\perp)e^\perp,
\]
we have that $e,-e\in \cB(\RR^d\setminus \hat{y}(\partial\Omega))$. 
Moreover, $F$ keeps the line in direction $e$ including those two points fixed, and
they are regular values for both $F$ and $y$. %, $y^{-1}(\{(-1,0)\})=\{(-1,0)\}$ and $y^{-1}(\{(2,0)\})=\{(0,0)\}$.
Hence, 
$\Deg(\hat{y};\Omega;e)=+\Deg(y;\Omega;e)=1$ and $\Deg(\hat{y};\Omega;-e)=-\Deg(y;\Omega;-e)=-1$,
because $\det \nabla F(-e)<0<\det \nabla F(e)$.
\end{ex}
\begin{rem}\label{rem:top3}
Starting with a domain with several holes, a similar construction as in Example~\ref{ex:top1} 
with a subsequent orientation-preserving deformation can also cause the images of holes to be stacked inside of each other. In fact, this way, 
with $\abs{n}$ holes for any given $n\in \ZZ$, we can get a deformation $y$ invertible on $\partial\Omega$,
such that $\Deg(y;\Omega;\cdot)$ attains the value $n$ somewhere. 
This also works in context of Example~\ref{ex:top2} if we have several small conical holes that all meet at the tip of the big outer cone.
\end{rem}
%\begin{rem}\label{rem:top4}
%It is also possible to construct examples satisfying $\det\nabla y\geq 0$ a.e.~in $\Omega$.
%The basic building block is again the map $y$ of Example~\ref{ex:top2}. Its determinant changes sign, but a negative sign occurs only in points that are mapped
%to $\cU$. Hence, composing $y$ with the parallel projection $S$ along $e^\perp$ onto $\bar\cB=\bar V_2\cup (-\bar{V}_1)$ which is globally Lipschitz on $y(\Omega)\setminus \bar\cB$,
%get $\tilde{y}:=S\circ y$ with the property that $\tilde{y}(\Omega)=V_2\cup (-V_1)$, $\det\nabla \tilde{y}\geq 0$ and $\det\nabla \tilde{y}(x)>0$ for every $x$ with $y(x)\in \cB$.
%%DOES NOT YET CONTRATICT THE THEOREM: DEG in 0,1. HOW TO COMBINE MORE THAN ONE HOLE?
%\end{rem}

\subsection*{Acknowledgements}

This research was supported by the Czech Science Foundation (GA \v{C}R) and the Austrian Science Fund (FWF) through the
bilateral grant 19-29646L (Large Strain Challenges in Materials Science), and through the associated 
MSMT-WTZ bilateral travel grant 8J19AT013.

%the grants GA18-03834S 

\bibliographystyle{plain}
\bibliography{../NLEbib}

\end{document}